\documentclass[12pt, reqno]{amsart} 

\usepackage{xcolor} 
\usepackage[colorlinks=true,linkcolor=blue,citecolor=red,backref=page]{hyperref} 
\usepackage{esint} 
\usepackage{amssymb} 
\usepackage{enumitem} 
\usepackage{mathrsfs} 

\setenumerate{label=(\alph*)} 
\setenumerate[2]{label=(\Alph*)} 

\newtheorem{theorem}{Theorem}[section]
\newtheorem{proposition}[theorem]{Proposition}
\newtheorem{lemma}[theorem]{Lemma}
\newtheorem{corollary}[theorem]{Corollary}

\theoremstyle{definition}
\newtheorem{definition}[theorem]{Definition}
\newtheorem{example}[theorem]{Example}
\theoremstyle{remark}
\newtheorem{remark}[theorem]{Remark}
\newtheorem{notationconvention}[theorem]{Notation Convention}

\numberwithin{equation}{section} 

\def\div{\operatorname{div}}

\def\Re{\operatorname{Re}}

\def\loc{\operatorname{loc}}
\def\inf{\operatornamewithlimits{inf\vphantom{p}}}

\newcommand{\norm}[1]{\left\lVert#1\right\rVert}
\newcommand\restr[2]{\ensuremath{\left.#1\right|_{#2}}}

\DeclareMathOperator*{\essinf}{ess\,inf}

\calclayout

\begin{document}

\title[The higher order regularity problem for elliptic systems]{The higher order regularity problem for elliptic systems with data in Banach function spaces}

\author{Juan Jos\'{e} Mar\'{i}n}
\address{Juan Jos\'{e} Mar\'{i}n
\\
Departamento de Matem{\'a}tica Aplicada y Estad{\'i}stica
\\
Universidad Polit{\'e}cnica de Cartagena
\\
Cartagena, Spain.} \email{juanjose.marin@upct.es}





\begin{abstract}
We characterize the well-posedness of the higher order regularity problem in the upper half-space with data in Sobolev Banach function spaces by proving its equivalence to natural weighted estimates for the Hardy–Littlewood maximal operator.
The generality of our framework allows for applications to Lebesgue spaces, rearrangement-invariant spaces such as Orlicz spaces, and variable exponent Lebesgue spaces, as well as their weighted counterparts, among others.
This is established for the family of second-order, homogeneous, elliptic, constant complex coefficient systems in ${\mathbb{R}}^n$ that admit a distinguished coefficient tensor, a natural condition that always holds in the scalar case and for the Lamé system of elasticity.
\end{abstract}

\maketitle
\tableofcontents

\section{Introduction and statement of main results}

First, fix two positive integers $M,n \in \mathbb{N}$ with $n\geq 2$ and consider the second-order, homogeneous, constant complex coefficient, $M\times M$ system $L$ in $\mathbb{R}^n$, defined as
\begin{equation}\label{eq:L-1}
Lu:=\left(a_{jk}^{\alpha\beta}\partial_j\partial_k u_\beta \right)_{1\leq\alpha,\beta\leq M},
\end{equation}
where $u=\left( u_\beta \right)_{1 \leq \beta \leq M}$ is a vector-valued function and
$a_{jk}^{\alpha\beta}\in\mathbb{C}$ for every $j,k\in\{1,...,n\}$ and $\alpha,\beta\in\{1,...,M\}$.
Here, and throughout this work, we make use of the summation convention over repeated indices.
Assume next that $L$ is (strongly) {\tt elliptic}, that is,
there exists a constant $c_0 \in (0,\infty)$ such that the following Legendre-Hadamard ellipticity condition holds: 
\begin{equation}\label{eq:L-2}
\begin{array}{c}
\Re \big[ a_{jk}^{\alpha\beta}\xi_j\xi_k\zeta_\alpha\zeta_\beta \big]
\geq c_0\left|\xi\right|^2\left|\zeta\right|^2\,\,\text{ for every}
\\[5pt]
\xi=\left(\xi_j\right)_{1\leq j\leq n}\in\mathbb{C}^n
\,\,\text{ and }\,\,\zeta=\left(\zeta_\alpha\right)_{1\leq\alpha\leq M}\in\mathbb{C}^M.
\end{array}
\end{equation}
Relevant examples include the Laplacian, as well as other more general scalar operators, such as those of the form ${\rm div}A\nabla$, where 
$A=(a_{jk})_{1\leq j,k\leq n}$ is an $n\times n$ matrix with constant complex coefficients satisfying the 
ellipticity condition
\begin{equation}\label{eq:intro-ajk}
\Re\big[ a_{jk}\xi_j\xi_k \big] \geq c_0\left|\xi\right|^2,\quad
\text{ for all }
\xi=\left(\xi_j\right)_{1\leq j\leq n}\in\mathbb{C}^n.
\end{equation}
Among non-scalar operators, a notable example is the Lamé system of elasticity (see Example \ref{ex:Lame}).

The central purpose of this work is to establish the well-posedness of the higher order regularity problem in the upper half-space $\mathbb{R}^n_+$ with data in Banach function spaces.
The key feature of this problem is that the boundary data belongs to a higher order Sobolev space and we demand nontangential control over higher order derivatives of the solution $u$. Specifically, we require estimates on $\mathcal{N}_{\kappa}(\nabla ^k u)$, where $\nabla^k u$ denotes the vector with components $(\partial^\alpha u)_{|\alpha|=k}$ and $\mathcal{N}_{\kappa}$ stands for the nontangential maximal function (see \eqref{eq:Nku} for a precise definition).

The uniqueness argument relies on a Fatou-type theorem presented in Section \ref{sec:Poisson}.
The proof of the existence of solutions is based on Poisson kernels. As shown by Agmon, Douglis, and Nirenberg in their pioneering works \cite{ADN59-1,ADN59-2}, any operator $L$ meeting the conditions \eqref{eq:L-1}–\eqref{eq:L-2} admits an associated Poisson kernel $P^L$ which serves as a natural generalization of the classical harmonic Poisson kernel (cf. Section \ref{sec:Poisson}).

However, as pointed out in \cite{MMMM-ho}, the higher order regularity problem is not expected to be well posed for arbitrary systems $L$ as above, even with boundary data in Lebesgue spaces. Consequently, we impose an additional algebraic assumption, namely that $L$ admits a distinguished coefficient tensor, $\mathfrak{A}_L^{\rm dis} \neq \varnothing$.
This condition is made precise in Section \ref{sec:tensors} and arises naturally in other settings related to boundary value problems. It is always satisfied in the scalar case $M=1$ and for systems such as the Lamé system of elasticity.

It should be emphasized that the results in this work are formulated in such a way that they include a broad range of (weighted and unweighted) Banach function spaces including Lebesgue spaces, rearrangement invariant spaces such as Orlicz spaces, and variable exponent Lebesgue spaces, among others (see Section \ref{sec:examples} for some examples).
For further details on Banach function spaces $\mathbb{X}_v$ and their associated spaces $\mathbb{X}'_v$, we refer the reader to Section \ref{sec:BFS}.
This general formulation requires a suitable definition for weighted Sobolev spaces. Specifically, if $\mathbb{M}$ denotes the space of measurable functions as in \eqref{eq:M}, $v$ and $w$ are weights on $\mathbb{R}^{n-1}$, and $\mathbb{X}_v$ is a Banach function space over $(\mathbb{R}^{n-1},v)$ then $(\mathbb{X}_{v})^w_k$ will stand for the $\mathbb{X}_v$-based Sobolev space of functions in $\mathbb{M}$ such that $fw \in \mathbb{X}_v$ and $(\partial^{\alpha} f)w \in \mathbb{X}_v$ for every multi-index $\alpha$ with $|\alpha| \leq k$ (see \eqref{eq:def-sobolev1}-\eqref{eq:def-sobolev2}).

We now present our main well-posedness result, whose proof can be found in Section~\ref{sec:proofs}. 

\begin{theorem}\label{thm:BVP-main}
Let $L$ be a second-order, homogeneous, elliptic, constant complex coefficient, $M\times M$ system in ${\mathbb{R}}^n$ as in \eqref{eq:L-1}-\eqref{eq:L-2} such that $\mathfrak{A}_L^{\rm dis} \neq \varnothing$. Let $v$ and $w$ be weights on $(\mathbb{R}^{n-1},\mathcal{L}^{n-1})$ such that $v \in L^1_{\loc}(\mathbb{R}^{n-1})$. Fix an aperture parameter $\kappa>0$, $\ell \in \mathbb{N}_0$ and a Banach function space $\mathbb{X}_v$ over $(\mathbb{R}^{n-1}, v)$. Let $\mathcal{M}$ denote the Hardy-Littlewood maximal function on $\mathbb{R}^{n-1}$ and let $\mathcal{M}'_{v}h:=\mathcal{M}(hv)/v$ for each $h \in\mathbb{M}$. Assume that there is $C_0\in(0,\infty)$ so that 
\begin{align}
\label{eq:M1} \norm{(\mathcal{M}h)\, w}_{\mathbb{X}_v} &\leq C_0\, \|hw\|_{\mathbb{X}_v}, \qquad \text{ for each }\,h \in\mathbb{M}, 
\\
\label{eq:M2} \|(\mathcal{M}'_{v} h)\, w^{-1} \|_{\mathbb{X}'_v} &\leq C_0\, \|hw^{-1}\|_{\mathbb{X}'_v}, \qquad \text{ for each }\,h \in\mathbb{M}.
\end{align}
Then, the $\ell$-th order regularity problem 
\begin{equation}\label{eq:BVP-main}
\left\lbrace
\begin{array}{l}
u\in\big[\mathscr{C}^{\infty}(\mathbb{R}^n_{+})\big]^M,
\\[4pt]
Lu=0\,\,\text{ in }\,\,\mathbb{R}^n_{+},
\\[4pt] 
(\mathcal{N}_{\kappa}(\nabla ^k u))\, w\in\mathbb{X}_v, \quad k\in\{0,1,...,\ell\},
\\[4pt]
\restr{u}{\partial\mathbb{R}^n_{+}}^{{}^{\kappa-{\rm n.t.}}}=f \in \big[(\mathbb{X}_v)^w_\ell \big]^M,
\end{array}
\right.
\end{equation}
is well-posed. More specifically, there is a unique solution and it is given by 
\begin{equation}\label{eq:BVP-main-sol}
u(x',t)=(P_t^L*f)(x'), \qquad(x',t)\in{\mathbb{R}}^n_{+},
\end{equation}
where $P^L$ denotes the Poisson kernel for $L$ in $\mathbb{R}^n_+$ from Theorem~\ref{thm:Poisson}. Furthermore, there is a constant $C\in[1,\infty)$ such that 
\begin{equation}\label{eq:BVP-main-estimate}
\sum_{k=0}^{\ell}\|(\mathcal{N}_{\kappa}(\nabla ^k u))\,w\|_{\mathbb{X}_v} \leq C\, \|f\|_{[ (\mathbb{X}_v)^w_\ell]^M}.
\end{equation}
\end{theorem}

The previous theorem is capable of establishing well-posedness results at every degree of regularity, thereby extending previous theories to encompass the full regularity spectrum.
Its broad formulation delivers new results even at the endpoint $\ell=0$, which corresponds to the Dirichlet problem.

\begin{remark}\label{rem:dual}
Since $(\mathbb{X}'_v)'=\mathbb{X}_v$ (cf. \eqref{eq:X''}), the hypothesis in Theorem \ref{thm:BVP-main} remain intact if we replace $\mathbb{X}_v$ by $\mathbb{X}'_v$ and $w$ by $w^{-1}v^{-1}$. Indeed, \eqref{eq:M1} becomes
\begin{equation}
\norm{(\mathcal{M}h)\, w^{-1}v^{-1}}_{\mathbb{X}'_v}
\leq C_0\, \|hw^{-1}v^{-1}\|_{\mathbb{X}'_v}, \qquad \text{ for each }\,h \in\mathbb{M}, 
\end{equation}
which is equivalent to \eqref{eq:M2} because $h\in\mathbb{M} \Leftrightarrow hv\in\mathbb{M}$ and,
conversely, \eqref{eq:M2} becomes \eqref{eq:M1}.
Therefore, under the hypothesis of Theorem \ref{thm:BVP-main}, the $\ell$-th regularity problem
\begin{equation}\label{eq:BVP-dual}
\left\lbrace
\begin{array}{l}
u\in\big[\mathscr{C}^{\infty}(\mathbb{R}^n_{+})\big]^M,
\\[4pt]
Lu=0\,\,\text{ in }\,\,\mathbb{R}^n_{+},
\\[4pt] 
(\mathcal{N}_{\kappa}(\nabla ^k u))\, w^{-1}v^{-1} \in\mathbb{X}'_v, \quad k\in\{0,1,...,\ell\},
\\[4pt]
\restr{u}{\partial\mathbb{R}^n_{+}}^{{}^{\kappa-{\rm n.t.}}}=f
\in \big[ (\mathbb{X}'_v)^{w^{-1}v^{-1}}_\ell \big]^M,
\end{array}
\right.
\end{equation}
also has a unique solution given by \eqref{eq:BVP-main-sol} and satisfying
\begin{equation}\label{eq:BVP-dual-estimate}
\sum_{k=0}^{\ell}\|(\mathcal{N}_{\kappa}(\nabla ^k u))\,w^{-1}v^{-1}\|_{\mathbb{X}'_v}
\leq C\, \|f\|_{[(\mathbb{X}'_v)^{w^{-1}v^{-1}}_\ell]^M} .
\end{equation}
\end{remark}

The following result shows that the results above are sharp in the sense that the well-posedness of the regularity problem \eqref{eq:BVP-main} and its dual counterpart \eqref{eq:BVP-dual} satisfying their naturally accompanying estimates \eqref{eq:BVP-main-estimate} and \eqref{eq:BVP-dual-estimate}
actually imply \eqref{eq:M1}-\eqref{eq:M2}.
In fact, we prove a stronger result assuming only well-posedness of the
Dirichlet problem (i.e., the case $\ell=0$) for the Laplacian. See Section \ref{sec:proofs} for the proof.

\begin{theorem}\label{thm:BVP-sharp}
Let $v$ and $w$ be weights on $(\mathbb{R}^{n-1},\mathcal{L}^{n-1})$ such that $v \in L^1_{\loc}(\mathbb{R}^{n-1})$. Fix an aperture parameter $\kappa>0$ and a Banach function space $\mathbb{X}_v$ over $(\mathbb{R}^{n-1}, v)$. Then, the following holds:
\begin{enumerate}
\item Assume that the Dirichlet problem:
\begin{equation}\label{eq:BVP-sharp1}
\left\lbrace
\begin{array}{l}
u\in\big[\mathscr{C}^{\infty}(\mathbb{R}^n_{+})\big]^M,
\\[4pt]
\Delta u=0\,\,\text{ in }\,\,\mathbb{R}^n_{+},
\\[4pt] 
(\mathcal{N}_{\kappa}u)\, w\in\mathbb{X}_v,
\\[4pt]
\restr{u}{\partial\mathbb{R}^n_{+}}^{{}^{\kappa-{\rm n.t.}}}=f \,
\text{ on }\mathbb{R}^{n-1},\quad fw\in\mathbb{X}_v,
\end{array}
\right.
\end{equation}
has a solution satisfying
\begin{equation}\label{eq:BVP-sharp1-estimates}
\qquad \|(\mathcal{N}_{\kappa}u)\,w\|_{\mathbb{X}_v}
\leq C\, \|fw\|_{\mathbb{X}_v}.
\end{equation}
Then, there is $C_0\in (0,\infty)$ such that 
\begin{equation} \label{eq:sharp1-M1}
\norm{(\mathcal{M}h)\, w}_{\mathbb{X}_v} \leq C_0\, \|hw\|_{\mathbb{X}_v}, \quad \text{ for each }\,h \in\mathbb{M}, 
\end{equation}
where $\mathcal{M}$ denotes the Hardy-Littlewood maximal function on $\mathbb{R}^{n-1}$.
\vspace*{0.1cm}
\item Assume that the Dirichlet problem:
\begin{equation}\label{eq:BVP-sharp2}
\left\lbrace
\begin{array}{l}
u\in\big[\mathscr{C}^{\infty}(\mathbb{R}^n_{+})\big]^M,
\\[4pt]
\Delta u=0\,\,\text{ in }\,\,\mathbb{R}^n_{+},
\\[4pt] 
(\mathcal{N}_{\kappa}u)\, w^{-1}v^{-1} \in\mathbb{X}'_v,
\\[4pt]
\restr{u}{\partial\mathbb{R}^n_{+}}^{{}^{\kappa-{\rm n.t.}}}=f \,
\text{ on }\mathbb{R}^{n-1},\quad fw^{-1}v^{-1}\in \mathbb{X}'_v,
\end{array}
\right.
\end{equation}
has a solution satisfying
\begin{equation}\label{eq:BVP-sharp2-estimates}
\|(\mathcal{N}_{\kappa}u)\,w^{-1}v^{-1}\|_{\mathbb{X}'_v}
\leq C\, \|fw^{-1}v^{-1}\|_{\mathbb{X}'_v}.
\end{equation}
Then, there is $C_0\in (0,\infty)$ such that 
\begin{equation} \label{eq:sharp2-M2}
\|(\mathcal{M}'_{v} h)\, w^{-1} \|_{\mathbb{X}'_v} \leq C_0\, \|hw^{-1}\|_{\mathbb{X}'_v}, \quad \text{ for each }\,h \in\mathbb{M}.
\end{equation}
where $\mathcal{M}$ denotes the Hardy-Littlewood maximal function on $\mathbb{R}^{n-1}$ and $\mathcal{M}'_{v}h:=\mathcal{M}(hv)/v$.
\vspace*{0.1cm}
\item The existence of a solution for each of the Dirichlet problems \eqref{eq:BVP-sharp1} and \eqref{eq:BVP-sharp2}, along with their corresponding estimates \eqref{eq:BVP-sharp1-estimates} and \eqref{eq:BVP-sharp2-estimates} is equivalent to conditions \eqref{eq:sharp1-M1} and \eqref{eq:sharp2-M2} together.
%
%
\end{enumerate}
\end{theorem}

As shown in the theorem above, the existence of a solution for the Laplacian Dirichlet problems \eqref{eq:BVP-sharp1} and \eqref{eq:BVP-sharp2} with their naturally accompanying estimates is enough to prove \eqref{eq:M1}-\eqref{eq:M2}. In turn, these bounds imply, by Theorem \ref{thm:BVP-main}, that the $\ell$-th order regularity problem \eqref{eq:BVP-main} is well-posed
for every $\ell \in \mathbb{N}_0$ and every system $L$ as in \eqref{eq:L-1}-\eqref{eq:L-2} with $\mathfrak{A}_L^{\rm dis}\neq \varnothing$. As a consequence, we obtain the following corollary.
\begin{corollary}\label{cor:intro}
Let $v$ and $w$ be weights on $(\mathbb{R}^{n-1},\mathcal{L}^{n-1})$ such that $v \in L^1_{\loc}(\mathbb{R}^{n-1})$. Fix an aperture parameter $\kappa>0$ and a Banach function space $\mathbb{X}_v$ over $(\mathbb{R}^{n-1}, v)$. Assume that each of the harmonic Dirichlet problems in
\eqref{eq:BVP-sharp1} and \eqref{eq:BVP-sharp2} has a solution satisfying its respective accompanying estimate \eqref{eq:BVP-sharp1-estimates} and \eqref{eq:BVP-sharp2-estimates}.
Then, for every second-order, homogeneous, elliptic, constant complex coefficient, $M\times M$ system $L$ in ${\mathbb{R}}^n$ as in \eqref{eq:L-1}-\eqref{eq:L-2} such that $\mathfrak{A}_L^{\rm dis} \neq \varnothing$ and every $\ell \in \mathbb{N}_0$, the $\ell$-th order regularity problem \eqref{eq:BVP-main} has a unique solution given by \eqref{eq:BVP-main-sol} and satisfying \eqref{eq:BVP-main-estimate}.
\end{corollary}

This section concludes with a brief discussion of some related results in the literature.
Specifically, we focus on boundary value problems in the upper half-space in which the boundary traces are taken nontangentially and the solutions are controlled using the nontangential maximal function.
In this vein, numerous contributions have addressed boundary value problems in Lebesgue spaces for the Laplacian (cf. \cite{ABR01,GR85,S70}) or, more generally, for variable coefficient scalar elliptic operators in divergence form, $L=\div A \nabla$, under different hypothesis on the coefficient matrix $A$ (cf. \cite{AAAHK11, HMM15}).
Regarding non-scalar operators, Fatou-type results and boundary value problems for
second-order, homogeneous, elliptic, constant complex coefficient systems have been extensively treated in \cite{C23-CMO, MMMMM19-FATOUPOISSON, MMM20-BMOW, MMMM-ho, MMMM16-Kothe, MMMM17-semigroup, MMMM19-BMO, MMMM19}.

Many of these works deal with the Dirichlet problem. In particular,
the Dirichlet problem in the upper half-space was treated in \cite{MMMM16-Kothe} within the context of Banach function spaces. In that article, the authors prove a general theorem from which they derive a result that covers the unweighted cases and a separate result, under the additional assumption that the function space is rearrangement invariant, to include weighted scenarios such as the one in Example \ref{ex:Lpw} (with $\ell=0$). The approach adopted here encompasses both the unweighted and weighted settings within a unified theorem and without requiring the rearrangement invariant assumption on $\mathbb{X}$.
In \cite{CMM22}, the Dirichlet problem for the Laplacian was analyzed in the context of Banach function spaces, making use of the same weighted formulation that we employ here.
Additionally, it has been shown in \cite{MMM21-failure} that the $L^p$-Dirichlet problem may fail to be well-posed if the ellipticity assumption in \eqref{eq:L-2} is replaced by a weak ellipticity condition.

The higher order regularity problem is treated in \cite{MMMM-ho} for $\mathbb{X}_v=L^p(\mathbb{R}^{n-1})$ in the unweighted case $v=w\equiv 1$. Moreover, in \cite{MMMM17-semigroup}, the authors establish a link between the infinitesimal generator of the Poisson semigroup and the solvability of the regularity problem in $L^p(\mathbb{R}^{n-1})$.

Similarly, well-posedness is proved for the regularity problem (of order one, i.e., $\ell=1$) in \cite[Theorem 8.19]{MMMMM22-USKT} using different techniques, for a class of domains that includes the upper half-space. While the Banach function spaces considered are of the form $\mathbb{X}_v$, the theorem does not cover certain cases such as the result for weighted Lebesgue spaces in Example \ref{ex:Lpw}.

\section{Preliminaries}\label{sec:prelim}

We start by introducing some notation and definitions.
As usual, $\mathbb{N}$ will denote the set of all strictly positive integers and 
$\mathbb{N}_0 := \mathbb{N} \cup \{ 0 \}$.
Next, consider the upper half-space
\begin{equation}
{\mathbb{R}}^n_{+} = \{ x = (x',x_n) \in \mathbb{R}^n
= \mathbb{R}^{n-1} \times \mathbb{R} : x_n > 0 \},
\end{equation}
with $n\in\mathbb{N}$ and $n\geq 2$.
Its boundary, $\partial{\mathbb{R}}^{n}_{+} = \mathbb{R}^{n-1} \times \{ 0 \}$, will be usually identified
with $\mathbb{R}^{n-1}$.
Next, fix an aperture parameter $\kappa>0$ and define, for each boundary point
$x'\in\partial{\mathbb{R}}^{n}_{+}$, the conical nontangential approach region
with vertex at $x'$ as
\begin{equation}
\Gamma_\kappa(x'):=\big\{y=(y',t)\in{\mathbb{R}}^{n}_{+}:\,|x'-y'|<\kappa\,t\big\}.
\end{equation}
For every function $u:{\mathbb{R}}^{n}_{+}\to{\mathbb{C}}^M$, we define
the nontangential maximal function of $u$ as
\begin{equation}\label{eq:Nku}
\big({\mathcal{N}}_\kappa u\big)(x'):=\sup\big\{|u(y)|:\,y\in\Gamma_\kappa(x')\big\},\qquad \text{ for }
x'\in\partial{\mathbb{R}}^{n}_{+}\equiv{\mathbb{R}}^{n-1}.
\end{equation}
Whenever meaningful, we also define the nontangential trace of $u$ as 
\begin{equation}\label{eq:trace}
u\big|^{{}^{\kappa-{\rm n.t.}}}_{\partial{\mathbb{R}}^{n}_{+}}(x')
:=\lim_{\Gamma_{\kappa}(x')\ni y\to (x',0)}u(y), \qquad
\quad\text{for }\,x'\in\partial{\mathbb{R}}^{n}_{+}\equiv{\mathbb{R}}^{n-1}.
\end{equation}

The Lebesgue measure on $\mathbb{R}^{n-1}$ will be denoted by $\mathcal{L}^{n-1}$ and $\mathbb{M}$ will stand for the set of all Lebesgue measurable functions which are finite almost everywhere, that is,
\begin{align}\label{eq:M}
\mathbb{M} = \{ f:\mathbb{R}^{n-1}\to [-\infty,\infty] & : \, 
f \text{ is Lebesgue measurable and } 
\nonumber \\ & \qquad
|f(x')|<\infty \, \text{ for }\, \mathcal{L}^{n-1} \text{-a.e. } x'\in \mathbb{R}^{n-1} \}.
\end{align}

Next, for every function $f$ in $\mathbb{R}^{n-1}$ whose components are locally integrable, the Hardy-Littlewood maximal operator on $\mathbb{R}^{n-1}$ 
is defined by 
\begin{equation}\label{eq:maximal}
\mathcal{M} f(x'):=\sup_{B \ni x'}\fint_{B}|f(y')|\,dy',\qquad 
\text{ for each }\,x'\in\mathbb{R}^{n-1},
\end{equation}
where the supremum is taken over all balls in $\mathbb{R}^{n-1}$ containing
$x'$, the
barred integral denotes mean average, and we take componentwise averages if
$f$ is $\mathbb{C}^M$-valued.


A {\tt weight} $w$ on $(\mathbb{R}^{n-1},\mathcal{L}^{n-1})$ is a measurable function satisfying $0<w(x')<\infty$ 
for $\mathcal{L}^{n-1}$-a.e. point $x'\in\mathbb{R}^{n-1}$. For $p\in(1,\infty)$, we say that $w$ belongs to the 
{\tt Muckenhoupt class} $A_p(\mathbb{R}^{n-1})$ if 
\begin{equation}\label{eq:def-Apweight}
[w]_{A_p}:=\sup_{B\subset\mathbb{R}^{n-1}}\left(\fint_{B}w(x')\,dx' \right)
\left(\fint_{B}w(x')^{1-p'}dx' \right)^{p-1}<\infty,
\end{equation}
where $p'$ is the conjugate exponent of $p$, i.e., $1/p+1/p'=1$, and the supremum runs over all 
balls $B$ in $\mathbb{R}^{n-1}$. Similarly, we say that $w\in A_1(\mathbb{R}^{n-1})$ if 
\begin{equation}\label{eq:def-A1weight}
[w]_{A_1}:=\sup_{B\subset\mathbb{R}^{n-1}}\left(\essinf_{x' \in B}w(x')\right)^{-1}
\left(\fint_{B}w(x')\,dx' \right)<\infty.
\end{equation}
In particular, 
\begin{equation}
\mathcal{M}w(x')\leq[w]_{A_1}w(x') \qquad \text{for }\mathcal{L}^{n-1}\text{-a.e. }x'\in\mathbb{R}^{n-1}.
\end{equation}
We also define 
\begin{equation}
A_{\infty}(\mathbb{R}^{n-1}) :=
\bigcup_{p \geq 1} A_p(\mathbb{R}^{n-1}).
\end{equation}
From \eqref{eq:def-Apweight}, one can easily check that
\begin{equation}\label{eq:Ap-p1}
\parbox{14cm}{\begin{center}
for $p\in(1,\infty)$, $w \in A_p(\mathbb{R}^{n-1})$ if and only if $w^{1-p'}\in A_{p'}(\mathbb{R}^{n-1})$
\end{center}}
\end{equation}
and Hölder's inequality yields
\begin{equation}\label{eq:Ap-p2}
\parbox{14cm}{\begin{center}
if $1\leq p<q\leq \infty$ then $A_p(\mathbb{R}^{n-1}) \subset A_q(\mathbb{R}^{n-1})$.
\end{center}}
\end{equation}

We will use the symbol $w$ for the weighted measure $w\,\mathcal{L}^{n-1}$ and
we will denote the 
weighted Lebesgue space
\begin{equation}
L^p(\mathbb{R}^{n-1},w):=L^p(\mathbb{R}^{n-1},w\,\mathcal{L}^{n-1}).
\end{equation}
As for unweighted Lebesgue spaces, we agree to write
\begin{equation}
L^p(\mathbb{R}^{n-1}):=L^p(\mathbb{R}^{n-1},\mathcal{L}^{n-1}).
\end{equation}
It is well-known that:
\begin{equation}\label{eq:Mbounded-Lpw}
\parbox{9.0cm}{\begin{center}
for any $p\in(1,\infty)$, ${\mathcal{M}}$ is bounded on $L^p(\mathbb{R}^{n-1},w)$ if and only if $w \in A_p(\mathbb{R}^{n-1}).$
\end{center}}
\end{equation}

\section{Banach function spaces}\label{sec:BFS}

In this section, we introduce Banach function spaces and outline several of their properties that will be instrumental later.

\begin{definition}\label{def:BFS}
Given a weight $v$ on $(\mathbb{R}^{n-1},\mathcal{L}^{n-1})$ and a mapping $\norm{\cdot}_{\mathbb{X}_v}:\mathbb{M} \to [0,\infty]$,
the set 
\begin{equation} 
\mathbb{X}_v = \{ f \in \mathbb{M}  : \norm{f}_{\mathbb{X}_v} < \infty \} 
\end{equation} 
is called {\tt Banach function space} over $(\mathbb{R}^{n-1},v)$ if $\norm{\cdot}_{\mathbb{X}_v}$ satisfies the following 
properties for all $f,g \in \mathbb{M}$:
\begin{enumerate}
\item $\norm{f}_{\mathbb{X}_v} = \norm{ \left| f \right| }_{\mathbb{X}_v}$ and $\norm{f}_{\mathbb{X}_v}=0$
if and only if $f=0$ $v$-a.e. on $\mathbb{R}^{n-1}$.
\vspace*{0.15cm}
\item $\norm{f+g}_{\mathbb{X}_v} \leq \norm{f}_{\mathbb{X}_v} + \norm{g}_{\mathbb{X}_v}$.
\vspace*{0.15cm}
\item $\norm{af}_{\mathbb{X}_v} = |a| \norm{f}_{\mathbb{X}_v}$ for all $a \in \mathbb{C}$.
\vspace*{0.15cm}
\item If $|f| \leq |g|$ $v$-a.e. on $\mathbb{R}^{n-1}$, then $\norm{f}_{\mathbb{X}_v} \leq \norm{g}_{\mathbb{X}_v}$.
\vspace*{0.15cm}
\item If $\{ f_i \}_i \subset \mathbb{M}$ is a sequence such that $|f_i|$ increases to 
$|f|$ $v$-a.e. on $\mathbb{R}^{n-1}$ as $i\to\infty$, then $\norm{f_i}_{\mathbb{X}_v}$ increases to $\norm{f}_{\mathbb{X}_v}$ as $i\to\infty$.
\vspace*{0.15cm}
\item If $E \subset \mathbb{R}^{n-1}$ is a measurable set of finite measure then its
characteristic function $\mathbf{1}_E$ satisfies
$\norm{\mathbf{1}_E}_{\mathbb{X}_v} < \infty$ and $\int_E |f(x')|\, dx' \leq C_E \norm{f}_{\mathbb{X}_v}$ for some $C_E \in (0,\infty)$.
\end{enumerate}
\end{definition}

It can be verified that if the mapping $\norm{\cdot}_{\mathbb{X}_v}$
satisfies the conditions stated above, then $(\mathbb{X}_v,\norm{\cdot}_{\mathbb{X}_v})$ is a complete normed space in accordance with the terminology.
Also, we wish to remark that in the definition above, as usual, we use the symbol $v$ for the weighted measure $v\,\mathcal{L}^{n-1}$. The case $v \equiv 1$ gives Banach function spaces $\mathbb{X}$ over $(\mathbb{R}^{n-1},\mathcal{L}^{n-1})$.

Given a Banach function space $\mathbb{X}_v$, we define its associate space, denoted by $\mathbb{X}'_v$, as the Banach function space defined through the norm
\begin{equation}\label{eq:associatedX}
\norm{f}_{ \mathbb{X}_v^{'}} := \sup \left\{ \int_{\mathbb{R}^{n-1}} |f(x')g(x')|\, v(x') \,dx' : 
g \in  \mathbb{X}_v, \enskip \norm{g}_{\mathbb{X}_v} \leq 1 \right\}.
\end{equation}
As a result, the following generalized H\"{o}lder inequality holds:
\begin{equation} \label{eq:XHolder}
\int_{\mathbb{R}^{n-1}} |f(x')g(x')| \, v(x') \, dx' \leq \norm{f}_{\mathbb{X}_v} \norm{g}_{\mathbb{X}_v^{'}},
\end{equation}
for every $f \in \mathbb{X}_v$ and $g \in \mathbb{X}_v^{'}$. Every Banach function space 
coincides with its second associated space (cf. \cite[Theorem 2.7, p.10]{BS88}), that is, 
\begin{equation}\label{eq:X''}
\mathbb{X}_v=(\mathbb{X}_v')'.
\end{equation}
As a consequence, 
\begin{equation} \label{eq:normXassociate}
\norm{f}_{ \mathbb{X}_v} = \sup \left\{ \int_{\mathbb{R}^{n-1}} |f(x')g(x')| \, v(x') \, dx' 
: g \in  \mathbb{X}'_v, \enskip \norm{g}_{\mathbb{X}'_v} \leq 1 \right\}.
\end{equation}

The following extrapolation theorem for Banach function spaces was proved in \cite{CMM22} in a more general context.
Remember that a family of extrapolation pairs $\mathcal{F}$ is just a collection of pairs $(f, g)$ of nonnegative measurable functions.

\begin{theorem}\label{thm:XExtrapolation}
Let $\mathcal{F}$ be a family of extrapolation pairs. Assume that $v$ and $w$ are weights on $(\mathbb{R}^{n-1}, \mathcal{L}^{n-1})$ such that $v \in L^1_{\loc}(\mathbb{R}^{n-1})$, and let $\mathbb{X}_v$ be a Banach function space over $(\mathbb{R}^{n-1}, v)$. Let $\mathcal{M}$ denote the Hardy-Littlewood maximal function on $\mathbb{R}^{n-1}$ and let $\mathcal{M}'_{v}h:=\mathcal{M}(hv)/v$ for each $h \in\mathbb{M}$. Fix $p_0\in (1,\infty)$, and assume that there are $C_1,C_2\in (0,\infty)$ such that 
\begin{align}
\|(\mathcal{M}h)\,w\|_{\mathbb{X}_v} &\leq C_1 \|hw\|_{\mathbb{X}_v}, \qquad \text{ for each }\,h \in \mathbb{M}, 
\\
\|(\mathcal{M}'_{v}h)\,w^{-1} \|_{\mathbb{X}'_v} &\leq C_2 \|hw^{-1}\|_{\mathbb{X}'_v}, \qquad\text{ for each }\,h \in\mathbb{M}.
\end{align}
If for every $\nu \in A_{p_0}(\mathbb{R}^{n-1})$, one has
\begin{equation}
\|f\|_{L^{p_0}(\mathbb{R}^{n-1},\nu)} \leq \Psi([\nu]_{A_{p_0, B}})\|g\|_{L^{p_0}(\mathbb{R}^{n-1},\nu)}, \qquad \text{ for each }(f,g) \in \mathcal{F}, 
\end{equation}  
where $\Psi: [1, \infty) \to [1, \infty)$ is a non-decreasing function, then 
\begin{equation}
\|fw\|_{\mathbb{X}_v} \leq C_0\, \|g\,w\|_{\mathbb{X}_v}, \qquad \text{ for each }(f,g) \in \mathcal{F},
\end{equation}
where $C_0=2^{3+\frac{2}{p'_0}}\,\Psi(2^{p_0}\,C_1^{p_0-1}\, C_2)$.
\end{theorem}

The next result was proved in \cite[(5.13)]{CMM22}
(see also \cite[Proposition 5.2.7]{MMM22-book2} and \cite[Proposition 8.2]{MMMMM22-USKT}).

\begin{proposition}\label{prop:inclusionX}
Let $v$ and $w$ be weights on $(\mathbb{R}^{n-1}, \mathcal{L}^{n-1})$ such that $v \in L^1_{\loc}(\mathbb{R}^{n-1})$, and let $\mathbb{X}_v$ be a Banach function space over $(\mathbb{R}^{n-1}, v)$. Let $\mathcal{M}$ denote the Hardy-Littlewood maximal function on $(\mathbb{R}^{n-1}, \mathcal{L}^{n-1})$ and let $\mathcal{M}'_{v}h:=\mathcal{M}(hv)/v$ for each $h \in\mathbb{M}$. Assume that there is $C_0\in (0,\infty)$ such that 
\begin{align}
\|(\mathcal{M}h)\,w\|_{\mathbb{X}_v} &\leq C_0 \|hw\|_{\mathbb{X}_v}, \qquad \text{ for each }\,h \in \mathbb{M}, 
\\
\|(\mathcal{M}'_{v}h)\,w^{-1} \|_{\mathbb{X}'_v} &\leq C_0 \|hw^{-1}\|_{\mathbb{X}'_v}, \qquad \text{ for each }\,h \in\mathbb{M}.
\end{align}
Then,
\begin{equation}\label{eq:hprop-2}
\int_{\mathbb{R}^{n-1}} \frac{|h(x')|}{1+|x'|^{n-1}}\,dx' < \infty, \qquad\text{for every $h \in \mathbb{M}$ with $\|hw\|_{\mathbb{X}_v}<\infty$}.
\end{equation}
\end{proposition}

We conclude this section by introducing higher order Sobolev spaces associated to a given Banach function space $\mathbb{X}_v$. 
Specifically, we will define a weighted version that is suitable for our purposes.
Given a weight $w$,
for each $k\in\mathbb{N}_0$, we define:
\begin{equation}\label{eq:def-sobolev1}
(\mathbb{X}_{v})^w_k = \{ f \in \mathbb{M} : fw \in \mathbb{X}_v \text{ and } (\partial^{\alpha} f)w \in \mathbb{X}_v
\text{ for each multi-index } \alpha \text{ with } |\alpha| \leq k \},
\end{equation}
which is a Banach space when equipped with the norm,
\begin{equation}\label{eq:def-sobolev2}
\displaystyle \norm{f}_{(\mathbb{X}_{v})^w_k} 
:= \norm{fw}_{\mathbb{X}_v} + \sum_{|\alpha|\leq k} \norm{(\partial^{\alpha} f)w}_{\mathbb{X}_v}.
\end{equation}

\section{Elliptic systems and distinguished coefficient tensors}\label{sec:tensors}

This section is devoted to compiling the essential properties of elliptic systems and the associated coefficient tensors. We begin by reviewing the most relevant aspects of fundamental solutions and, after that, we introduce the class of distinguished coefficient tensors along with some examples. We remark that much of the discussion presented here remains valid in the broader setting of weakly elliptic systems.

The following result gathers some features of the fundamental solution
associated to a given system $L$ as in \eqref{eq:L-1}-\eqref{eq:L-2}. For a proof, see \cite[Theorem 11.1]{M18-book}.

\begin{theorem}\label{thm:E}
Let $L$ be a second-order, homogeneous, elliptic, constant complex coefficient, $M\times M$ system in ${\mathbb{R}}^n$ as in \eqref{eq:L-1}-\eqref{eq:L-2}.
Then, there exists an $M\times M$ matrix-valued 
function $E=\big(E_{\alpha\beta}\big)_{1\leq \alpha,\beta\leq M}$, 
which we will call {\tt fundamental solution} associated with $L$,
such that the following properties hold.
\begin{enumerate}
\item For each $\alpha,\beta\in\{1,\dots,M\}$ one has 
$E_{\alpha\beta}\in{\mathscr{C}}^{\infty}\big(\mathbb{R}^{n}\setminus\{0\}\big)$ 
and $E_{\alpha\beta}(x)=E_{\alpha\beta}(-x)$ for every $x\in{\mathbb{R}}^n\setminus\{0\}$.

\vskip 0.08in
\item For every fixed point $y\in{\mathbb{R}}^n$ one has $L\big[E(\cdot-y)\big]=\delta_yI_{M\times M}$
in the sense of distributions in ${\mathbb{R}}^n$, where $I_{M\times M}$ is the $M\times M$ identity matrix 
and $\delta_y$ denotes the Dirac distribution in ${\mathbb{R}}^n$ with mass at $y$.
That is, using the standard Kronecker delta notation, 
\begin{equation}\label{eq:E-1}
a_{jk}^{\alpha\beta}\partial_{x_j}\partial_{x_k}\big[E_{\beta\gamma}(x-y)\big]=\delta_{\alpha\gamma}\delta_y(x),
\qquad \text{ for each }x\in{\mathbb{R}}^n,
\end{equation} 
in the sense of distributions, for every $\alpha,\gamma\in\{1,\dots,M\}$.

\vskip 0.08in
\item For each multi-index $\alpha\in{\mathbb{N}}_0^{n}$ with $n+|\alpha|>2$, the function 
$\partial^\alpha E$ is positive homogeneous of degree $2-n-|\alpha|$ and there exists a constant $C_\alpha\in(0,\infty)$ 
such that
\begin{equation}\label{E-2}
\big|(\partial^\alpha E)(x)\big|\leq C_\alpha|x|^{2-n-|\alpha|}\,\,\text{ for every }\,\,x\in{\mathbb{R}}^n\setminus\{0\}.
\end{equation}
Finally, if $n=2$, there is $C\in(0,\infty)$ such that 
$|E(x)|\leq C\big(1+\big|\ln|x|\big|\big)$ for every $x\in{\mathbb{R}}^2\setminus\{0\}$.
\end{enumerate}
\end{theorem}

Consider $A:=\big(a^{\alpha\beta}_{jk}\big)$, with $a^{\alpha\beta}_{jk} \in \mathbb{C}$ for every $j,k\in\{1,...,n\}$ and $\alpha,\beta\in\{1,...,M\}$,
and define the elliptic operator $L$ as in \eqref{eq:L-1}-\eqref{eq:L-2} for this choice of $A$.
We refer to $A$ as a {\tt coefficient} {\tt tensor} of $L$ and for convenience we will write $L=L_A$.
Next, consider two different coefficient tensors
\begin{equation}\label{eq:intro-AandtildeA}
A:=\big(a^{\alpha\beta}_{jk}\big)_{\substack{1\leq j,k\leq n\\1\leq\alpha,\beta\leq M}},
\qquad
\widetilde{A}:=\big(\widetilde{a}^{\,\alpha\beta}_{jk}\big)
_{\substack{1\leq j,k\leq n\\1\leq\alpha,\beta\leq M}},
\end{equation}
and assume that $a^{\alpha\beta}_{jk}-\widetilde{a}^{\,\alpha\beta}_{jk}$
is antisymmetric in the indices $j,k$ for each $\alpha,\beta\in\{1,\dots,M\}$.
Then, one can check that the systems generated by $A$ and $\widetilde{A}$ using \eqref{eq:L-1} turn out to be the same,
that is, $L_A=L_{\widetilde{A}}$. This demonstrates that a single system $L$ may be written using multiple coefficient tensors in \eqref{eq:L-1}.
Accordingly, it is natural to introduce the collection of all coefficient tensors associated with a fixed system $L$,
\begin{equation}\label{eq:UL}
\mathfrak{A}_L := \Big\{ A = \big( a_{jk}^{\alpha\beta} \big)
_{\substack{1\leq j,k\leq n\\1\leq\alpha,\beta\leq M}} : L = L_A \Big\}.
\end{equation}

We next record a result established in \cite[Proposition 1.2.1]{MMM23-book5}, which describes a number of equivalent properties for coefficient tensors.

\begin{proposition}\label{thm:distingcond}
Let $L$ be a second-order, homogeneous, elliptic, constant complex coefficient, $M\times M$ system in ${\mathbb{R}}^n$ as in \eqref{eq:L-1}-\eqref{eq:L-2}. For each 
$\xi\in{\mathbb{R}}^n\setminus\{0\}$, define the matrix-valued function $\mathcal{E}(\xi)$ as 
\begin{equation}
\big(\mathcal{E}_{\gamma\beta}(\xi)\big)_{1\leq\gamma,\beta\leq M}:=\left[
\left(-a_{jk}^{\alpha\beta}\xi_j\xi_k\right)_{1 \leq \alpha,\beta \leq M}
\right]^{-1}\in\mathbb{C}^{M\times M}.
\end{equation}
Let $E=\big(E_{\alpha\beta})_{1\leq\alpha,\beta\leq M}$ be the fundamental solution associated 
with $L$ as in Theorem~\ref{thm:E}. 
Then for each coefficient tensor 
$A=\big(a^{\alpha\beta}_{jk}\big)_{\substack{1\leq j,k\leq n \\ 1\leq\alpha,\beta\leq M}}\in{\mathfrak{A}}_L$ 
the following conditions are equivalent:
\begin{enumerate}
\item[(a)] For each $k,k'\in\{1,\dots,n\}$ and each $\alpha,\gamma\in\{1,\dots,M\}$ one has
\begin{equation}\label{eq:dis-a}
\big(x_{k'}a_{jk}^{\beta\alpha}-x_k a_{jk'}^{\beta\alpha}\big)(\partial_j E_{\gamma\beta})(x)=0
\quad\text{ for each }\,\,x\in\mathbb{R}^n \setminus\{0\}.
\end{equation}
\vspace*{0.1cm}
\item[(b)] For every $k,k'\in\{1,\dots,n\}$ and every $\alpha,\gamma\in\{1,\dots,M\}$,  one has
\begin{equation}\label{eq:dis-b}
\Big[a^{\beta\alpha}_{jk}\partial_{\xi_{k'}}-a^{\beta\alpha}_{jk'}\partial_{\xi_{k}}\Big]
\big[\xi_j{\mathcal{E}}_{\gamma\beta}(\xi)\big]=0.
\end{equation}
in the sense of tempered distributions in ${\mathbb{R}}^n$.
\vspace*{0.1cm}
\item[(c)] There exists a matrix-valued function 
\begin{equation}\label{eq:dis-c}
k=\big\{k_{\gamma\alpha}\big\}_{1\leq\gamma,\alpha\leq M}:{\mathbb{R}}^{n}\setminus\{0\}
\longrightarrow{\mathbb{C}}^{M\times M}
\end{equation}
with the property that for every $\gamma,\alpha\in\{1,\dots,M\}$ and $k\in\{1,\dots,n\}$, one has
\begin{equation}\label{eq:dis-c2}
a^{\beta\alpha}_{jk}(\partial_j E_{\gamma\beta})(x)=x_k k_{\gamma\alpha}(x)
\quad\text{ for each }\,\,x\in{\mathbb{R}}^n\setminus\{0\}.
\end{equation}
\end{enumerate}
\end{proposition}

\begin{definition}\label{def:dis}
Let $L$ be a second-order, homogeneous, elliptic, constant complex coefficient, $M\times M$ system in ${\mathbb{R}}^n$ as in \eqref{eq:L-1}-\eqref{eq:L-2}. We say that 
\begin{equation}\label{eq:dis-def1}
A=\big(a_{jk}^{\alpha\beta}\big)_{\substack{1\leq j,k\leq n\\ 1\leq\alpha,\beta\leq M}}\in{\mathfrak{A}}_L
\end{equation}
is a {\tt distinguished} {\tt coefficient} {\tt tensor}\index{distinguished coefficient tensor} 
for the system $L$ provided any of the equivalent conditions 
{\it (a)}-{\it (c)} in Proposition~\ref{thm:distingcond} holds. Also, denote by ${\mathfrak{A}}_L^{\rm dis}$ 
the family of such distinguished coefficient tensors for $L$, that is,
\begin{align}\label{eq:dis-def2}
{\mathfrak{A}}_L^{\rm dis}:=\Big\{A= &\big(a_{jk}^{\alpha\beta}\big)
_{\substack{1\leq j,k\leq n\\ 1\leq\alpha,\beta\leq M}}\in{\mathfrak{A}}_L:\,
\text{condition \eqref{eq:dis-a}}
\\[0pt]
& \text{holds for each }\,\,k,k'\in\{1,\dots,n\}\,\,\text{ and }\,\,\alpha,\gamma\in\{1,\dots,M\}\Big\}.
\nonumber
\end{align}
\end{definition}

Distinguished coefficient tensors appear naturally in the analysis of layer potentials. Specifically, in appropiate settings, every double layer potential admits a rich Calderón\-–Zygmund theory but requires demanding that the coefficient tensor is distinguished to obtain other relevant properties (cf. \cite{HMT10, MMMMM22-USKT, MMM23-book5}).
In the present framework, the condition ${\mathfrak{A}}_L^{\rm dis} \neq \varnothing$ also arises and it is actually a hypothesis in Theorem \ref{thm:BVP-main}. As a consequence, the following remark is relevant to ensure its applicability in the scalar case.

\begin{remark}\label{rem:scalardis}
%
Assume $L$ is a second-order, homogeneous, elliptic, constant complex coefficient, scalar differential operator 
in ${\mathbb{R}}^n$ (i.e., as in \eqref{eq:L-1}-\eqref{eq:L-2} with $M=1$).
In other words, suppose that $L=a_{jk}\partial_j\partial_k$ with $a_{jk}\in{\mathbb{C}}$ for $j,k\in\{1,\dots,n\}$ and
there is a constant $c_0\in(0,\infty)$ such that
\begin{equation}\label{eq:scalardis-1}
\Re\big[ a_{jk}\xi_j\xi_k \big] \geq c_0\left|\xi\right|^2,\quad
\text{ for all }
\xi=\left(\xi_j\right)_{1\leq j\leq n}\in\mathbb{C}^n.
\end{equation}
Define the coefficient matrices $A:=(a_{jk})_{1\leq j,k\leq n}\in{\mathbb{C}}^{n\times n}$ and
\begin{equation}\label{eq:scalardis-2} 
{\rm sym}(A):=\frac{A+A^{\top}}{2},
\end{equation}
and note that $L=L_A=L_{{\rm sym}(A)}$, that is,
${\rm sym}(A)$ can also be used to represent $L$ as in \eqref{eq:L-1}.
It turns out that ${\rm sym}(A)$ is actually the only distinguished coefficient tensor associated with $L$ (cf. \cite[Proposition 1.6.11]{MMM23-book5} for a proof). In particular,
${\rm sym}(A)\in{\mathfrak{A}}_L^{\rm dis}$ and hence,
\begin{equation}\label{eq:scalardis-6}
\parbox{11.00cm}{\begin{center}
$\mathfrak{A}_L^{\rm dis}\neq\varnothing$ for every second-order, homogeneous, elliptic, 
constant complex coefficient, scalar operator $L$ in ${\mathbb{R}}^n$.
\end{center}} 
\end{equation}
\end{remark}

We conclude this section by analyzing the Lamé system of elasticity as an example of an elliptic system which admits a distinguished coefficient tensor.

\begin{example}\label{ex:Lame}
Consider the complex form of the Lam\'e system, which arises in the theory of linear elasticity, defined in terms of two parameters $\mu,\lambda\in{\mathbb{C}}$ known as the Lamé moduli, as
\begin{equation}\label{eq:Lame-def}
L_{\mu,\lambda}:=\mu\Delta+(\mu+\lambda)\nabla\div,
\end{equation}
acting on vector-valued functions $u=(u_\beta)_{1\leq\beta\leq n}$ defined in ${\mathbb{R}}^n$ and
with the Laplacian applied componentwise. One can check
(cf., e.g., \cite{MMMM19}) that 
\begin{align}\label{eq:Lame-elliptic}
L_{\mu,\lambda}\text{ satisfies the ellipticity condition \eqref{eq:L-2}}\Longleftrightarrow 
\enskip \Re\mu>0,\enskip \Re (2\mu+\lambda)>0.
\end{align}
The Lamé system $L_{\mu,\lambda}$ can be represented in the form \eqref{eq:L-1}, with $M:=n$, 
using different coefficient tensors, including those of the one-parameter family 
\begin{equation}\label{eq:Lame-coef}
a_{jk}^{\alpha\beta}(\theta):=\mu\delta_{jk}\delta_{\alpha\beta}
+(\mu+\lambda-\theta)\delta_{j\alpha}\delta_{k\beta}+\theta\delta_{j\beta}\delta_{k\alpha},
\quad 1\leq j,k,\alpha,\beta\leq n.
\end{equation}
In other words, for each vector-valued function
$u=(u_\beta)_{1\leq\beta\leq n}$
and each parameter $\theta\in{\mathbb{C}}$, one has
\begin{equation}\label{eq:Lame-L}
L_{\mu,\lambda}u=\Big(a_{jk}^{\alpha\beta}(\theta)\partial_j\partial_k u_\beta\Big)_{1\leq\alpha\leq n}.
\end{equation}
It turns out that for any
$\mu,\lambda,\theta\in{\mathbb{C}}$ we have (cf. \cite{HMT10} or \cite[pp. 56-57]{MMM23-book5} for details),
\begin{equation}\label{eq:Lame-dis}
A(\theta)\in\mathfrak{A}_{L_{\mu,\lambda}}^{\rm dis}\enskip\,\Longleftrightarrow
\enskip\,3\mu+\lambda\neq0 \enskip \text{and}
\enskip\,\theta=\frac{\mu(\mu+\lambda)}{3\mu+\lambda}.
\end{equation}
In conclusion, \eqref{eq:Lame-elliptic} and \eqref{eq:Lame-dis} shows that the Lamé system of elasticity falls within the scope of this paper. Specifically,
\begin{equation}
\parbox{10.4cm}{\begin{center}
whenever $\Re \mu >0,\enskip \Re (2\mu+\lambda)>0$ and $3\mu+\lambda\neq0$,
$L_{\mu,\lambda}$ is a system as in \eqref{eq:L-1}-\eqref{eq:L-2} and satisfying $\mathfrak{A}_{L_{\mu,\lambda}}^{\rm dis} \neq \varnothing$.
\end{center}
}\end{equation}
\end{example}

\section{Poisson kernels and Fatou-type results}\label{sec:Poisson}

In this section we present the concept of Poisson kernel in ${\mathbb{R}}^n_{+}$
for a given system $L$ as in \eqref{eq:L-1}-\eqref{eq:L-2}. This generalizes
the classical harmonic Poisson kernel given by
\begin{equation}\label{eq:Poisson-harmonic}
P^{\Delta}(x') := \frac{2}{\omega_{n-1}}
\frac{1}{(1 + |x'|^2)^{n/2}}, \qquad x' \in \mathbb{R}^{n-1},
\end{equation}
where $\omega_{n-1}$ denotes the area of the unit sphere in $\mathbb{R}^n$.

\begin{definition}\label{def:Poisson}
Let $L$ be a second-order, homogeneous, elliptic, constant complex coefficient, $M\times M$ 
system in ${\mathbb{R}}^n$ as in \eqref{eq:L-1}-\eqref{eq:L-2}.
A {\tt Poisson kernel} for $L$ in $\mathbb{R}^{n}_{+}$ is a 
matrix-valued function 
\begin{equation}\label{eq:Poisson-def1}
P^L=\big(P^L_{\alpha\beta}\big)_{1\leq\alpha,\beta\leq M}:\mathbb{R}^{n-1}\longrightarrow\mathbb{C}^{M\times M}
\end{equation}
such that the following conditions hold:
\begin{list}{$(\theenumi)$}{\usecounter{enumi}\leftmargin=.8cm
\labelwidth=.8cm\itemsep=0.2cm\topsep=.1cm
\renewcommand{\theenumi}{\alph{enumi}}}
\item There exists $C\in(0,\infty)$ such that
\begin{equation}
|P^L(x')|\leq\frac{C}{(1+|x'|^2)^{\frac{n}2}} \quad \text{ for every } x'\in\mathbb{R}^{n-1}.
\end{equation}
\item The function $P^L$ is Lebesgue measurable and
\begin{equation}
\int_{\mathbb{R}^{n-1}}P^L(x')\,dx'=I_{M\times M},
\end{equation}
where $I_{M\times M}$ denotes the $M\times M$ identity matrix.
\item For each $x'\in\mathbb{R}^{n-1}$ and $t\in(0,\infty)$, define
\begin{equation}
K^L(x',t):=P^L_t(x'):=t^{1-n}P^L(x'/t).
\end{equation}
Then, the function
$K^L=\big(K^L_{\alpha\beta}\big)_{1\leq\alpha,\beta\leq M}$
satisfies
\begin{equation}\label{eq:Poisson-def2}
LK^L_{\cdot\beta}=0 \quad \text{ in } \,\mathbb{R}^{n}_{+}
\,\,\text{ for each }\,\,\beta\in\{1,\dots,M\},
\end{equation}
in the sense of distributions, where $K^L_{\cdot\beta}:=\big(K^L_{\alpha\beta}\big)_{1\leq\alpha\leq M}$.
\end{list}
\end{definition}

Poisson kernels for elliptic operators have received considerable attention as a tool for addressing boundary value problems in a half-space (see \cite{ADN59-1,ADN59-2,KMR01,S63,S64,S66}). 
The following result, derived as a corollary from a more general theory developed by S. Agmon, A. Douglis, and L. Nirenberg in \cite{ADN59-2}, ensures the existence of a Poisson kernel in our setting.

\begin{theorem}\label{thm:Poisson}
Let $L$ be a second-order, homogeneous, elliptic, constant complex coefficient, $M\times M$ 
system in ${\mathbb{R}}^n$ as in \eqref{eq:L-1}-\eqref{eq:L-2}. Then, there exists
a Poisson kernel $P^L$ associated with $L$ in the sense of Definition~\ref{def:Poisson}, which has the additional
property that the function defined as
\begin{equation}\label{eq:Poisson-1}
K^L(x',t):=P^L_t(x')\qquad\text{for each}\,\,(x',t)\in{\mathbb{R}}^n_{+},
\end{equation}
satisfies $K^L\in\big[{\mathscr{C}}^\infty\big(\overline{{\mathbb{R}}^n_{+}}
\setminus B(0,\varepsilon)\big)\big]^{M\times M}$ for every $\varepsilon>0$. 
Moreover, for each multi-index $\alpha\in{\mathbb{N}}_0^n$
there exists $C_\alpha\in(0,\infty)$ such that
\begin{equation}\label{eq:Poisson-2}
\big|(\partial^\alpha K^L)(x)\big|\leq C_\alpha\,|x|^{1-n-|\alpha|},\,\,\,\text{ for each }\,\,
x\in{\overline{{\mathbb{R}}^n_{+}}}\setminus\{0\}.
\end{equation}
\end{theorem}

Throughout this work, the convolution of a matrix-valued function with a vector-valued function is interpreted naturally in terms of matrix-vector multiplication. The next theorem, which appears in \cite[Theorem 3.1]{MMMM16-Kothe} is a fundamental tool to establish the existence of solutions for boundary value problems.

\begin{proposition}\label{prop:existence}
Let $L$ be a second-order, homogeneous, elliptic, constant complex coefficient, $M\times M$ 
system in ${\mathbb{R}}^n$ as in \eqref{eq:L-1}-\eqref{eq:L-2} and let $P^L$ be 
the Poisson kernel associated with $L$ in $\mathbb{R}^{n}_{+}$ from Theorem \ref{thm:Poisson}.
Fix an aperture parameter $\kappa>0$. Given a function 
\begin{equation}\label{eq:existence-f}
f\in\Big[L^1\Big(\mathbb{R}^{n-1}\,,\,\frac{dx'}{1+|x'|^n}\Big)\Big]^M,
\end{equation}
set
\begin{equation}\label{eq:existence-u}
u(x',t):=(P^L_t\ast f)(x'),\qquad\text{ for all }\,(x',t)\in{\mathbb{R}}^n_{+}.
\end{equation}
Then $u$ is meaningfully defined via an absolutely convergent integral,
\begin{equation}\label{eq:existence-u2}
u\in\big[\mathscr{C}^\infty(\mathbb{R}^n_{+})\big]^{M},\quad
Lu=0\,\text{ in }\,\mathbb{R}^{n}_{+},\,\,\, \text{ and } \,\,\,
u\big|_{\partial\mathbb{R}^{n}_{+}}^{{}^{\kappa-{\rm n.t.}}}=f
\,\,\text{ ${\mathcal{L}}^{n-1}$-a.e. on }\,\mathbb{R}^{n-1},
\end{equation}
where the last identity is valid in the set of Lebesgue points of $f$. Moreover,
there exists a constant $C\in(0,\infty)$ depending on $n$, $\kappa$, and $L$ such that
\begin{equation}\label{eq:existence-NuMf}
\big(\mathcal{N}_{\kappa}u\big)(x')\leq C\big(\mathcal{M}f\big)(x'),\qquad\text{ for each }\,x'\in\mathbb{R}^{n-1}.
\end{equation}
\end{proposition}

We finish this section with the following Fatou-type result from \cite[Theorem 1.1]{MMMMM19-FATOUPOISSON},
which will be used to establish the uniqueness of solutions for boundary value problems.

\begin{proposition}\label{prop:uniq}
Let $L$ be a second-order, homogeneous, elliptic, constant complex coefficient, $M\times M$ 
system in ${\mathbb{R}}^n$ as in \eqref{eq:L-1}-\eqref{eq:L-2}.
Assume that $u\in\big[{\mathscr{C}}^{\infty}(\mathbb{R}^n_{+})\big]^M$ satisfies $Lu=0$ in $\mathbb{R}^n_{+}$ and 
\begin{equation}\label{eq:uniq-hyp}
\int_{\mathbb{R}^{n-1}}\frac{({\mathcal{N}}_{\kappa}u)(x')}{1+|x'|^{n-1}}\,dx'<\infty.
\end{equation}
Then,
\begin{align}
\begin{cases}
u\big|^{{}^{\kappa-{\rm n.t.}}}_{\partial{\mathbb{R}}^n_{+}}\ \text{exists } 
\mathcal{L}^{n-1}\text{-a.e. on}\ \mathbb{R}^{n-1},
\\[12pt]
\displaystyle u\big|^{{}^{\kappa-{\rm n.t.}}}_{\partial{\mathbb{R}}^n_{+}} \in \Big[ L^1\Big(\mathbb{R}^{n-1},\frac{dx'}{1+|x'|^{n-1}}\Big) \Big]^M,
\\[14pt]
u(x',t)=\Big(P^L_t\ast\big(u\big|^{{}^{\kappa-{\rm n.t.}}}_{\partial \mathbb{R}^n_{+}}\big)\Big)(x')\ \text{for each}\ (x',t)\in \mathbb{R}^n_{+},
\end{cases}
\end{align}
where $P^L$ denotes the Poisson kernel for $L$ in $\mathbb{R}^n_+$ from Theorem~\ref{thm:Poisson}.
In particular, if $u\big|^{{}^{\kappa-{\rm n.t.}}}_{\partial \mathbb{R}^n_{+}}=0$ $\mathcal{L}^{n-1}$-a.e. on $\mathbb{R}^{n-1}$, then $u\equiv 0$.
\end{proposition}

\section{Proofs of Theorems \ref{thm:BVP-main} and \ref{thm:BVP-sharp}}\label{sec:proofs}

The purpose of this section is to prove the main results stated in the introduction. The general scheme of the proof of Theorem \ref{thm:BVP-main} proceeds along the original approach in \cite{MMMM-ho} for $L^p$ spaces, while the proof of Theorem \ref{thm:BVP-sharp} relies on some ideas in \cite{MMMM16-Kothe}.

We need to develop several preliminary tools first, beginning with the following definition of a class of auxiliary kernels.

\begin{definition}\label{def:Q-T}
Let $L$ be a second-order, homogeneous, elliptic, constant complex coefficient, $M\times M$ 
system in ${\mathbb{R}}^n$ as in \eqref{eq:L-1}-\eqref{eq:L-2}. Then, for every
index $j \in \{1,...,n\}$ we define the matrix-valued function
\begin{equation}
Q^{(j)}(x') := 
\left( Q^{(j)}_{\alpha\beta}(x') \right)_{1 \leq \alpha,\beta \leq M} :=
\left( (\partial_j E_{\alpha\beta}(x',1) \right)_{1 \leq \alpha,\beta \leq M},
\quad \text{ for each } x' \in \mathbb{R}^{n-1},
\end{equation}
where $E$ is the fundamental solution for $L$ (cf. Theorem \ref{thm:E}).
\end{definition}

We now proceed to prove an estimate for the nontangential maximal function of convolutions with the previously defined kernels, in the setting of Banach function spaces.

\begin{lemma}\label{lem:NT}
Let $L$ be a second-order, homogeneous, elliptic, constant complex coefficient, $M\times M$ 
system in ${\mathbb{R}}^n$ as in \eqref{eq:L-1}-\eqref{eq:L-2}. Fix $\alpha,\beta\in\{1,...,M\}$ and $j \in \{1,...,n\}$.
Let $v$ and $w$ be weights on $(\mathbb{R}^{n-1},\mathcal{L}^{n-1})$ such that $v \in L^1_{\loc}(\mathbb{R}^{n-1})$. Fix an aperture parameter $\kappa>0$ and a Banach function space $\mathbb{X}_v$ over $(\mathbb{R}^{n-1}, v)$. Let $\mathcal{M}$ denote the Hardy-Littlewood maximal function on $\mathbb{R}^{n-1}$ and let $\mathcal{M}'_{v}h:=\mathcal{M}(hv)/v$ for each $h \in\mathbb{M}$. Assume that there is $C_0\in(0,\infty)$ so that 
\begin{align}
\label{eq:NT-M1} \norm{(\mathcal{M}h)\, w}_{\mathbb{X}_v} &\leq C_0\, \|hw\|_{\mathbb{X}_v}, \qquad \text{ for each }\,h \in\mathbb{M}, 
\\
\label{eq:NT-M2} \|(\mathcal{M}'_{v} h)\, w^{-1} \|_{\mathbb{X}'_v} &\leq C_0\, \|hw^{-1}\|_{\mathbb{X}'_v}, \qquad \text{ for each }\,h \in\mathbb{M}.
\end{align}
Then, there exists a constant $C\in(0,\infty)$ such that
\begin{equation}\label{eq:NT}
\|(\mathcal{N}_{\kappa}(( Q_{\alpha\beta}^{(j)})_t \ast f)\,w\|_{\mathbb{X}_v} \leq C\, \|fw\|_{\mathbb{X}_v},
\qquad \text{ whenever } \,fw \in\mathbb{X}_v.
\end{equation}
\end{lemma}
\begin{proof}
Define
\begin{equation}
T^{(j)}_{\alpha\beta}f(x',t) :=
\left[ \left( Q_{\alpha\beta}^{(j)} \right)_t \ast f \right](x') =
\int_{\mathbb{R}_{n-1}} (\partial_j E_{\alpha\beta})(x'-y',t)\,f(y')\,dy',
\qquad (x',t) \in \mathbb{R}^{n}.
\end{equation}
From Theorem \ref{thm:E} we have that the kernel of this operator,
$\partial_j E_{\alpha\beta} \in \mathscr{C}^{\infty}(\mathbb{R}^n \setminus \{0\})$, is odd and positive homogeneous of degree
$1-n$, that is,
\begin{equation}
\partial_j E_{\alpha\beta}(-x)=-\partial_j E_{\alpha\beta}(x),
\quad \partial_j E_{\alpha\beta}(\lambda x)=\lambda^{1-n}\partial_j E_{\alpha\beta}(x)
\end{equation}
for every $x \in \mathbb{R}^n \setminus \{0\}$ and $\lambda \in (0,\infty)$.
As a consequence of this and standard estimates for Muckenhoupt weights
(cf., e.g., \cite[Proposition 3.4]{MMMMM22-USKT}), for every $p \in (1,\infty)$ and $\nu \in A_{p}(\mathbb{R}^{n-1})$ there is $C\in(0,\infty)$ such that
\begin{equation}\label{eq:NT-Lpw}
\|\mathcal{N}_{\kappa}(T^{(j)}_{\alpha\beta}f)\|_{L^p(\mathbb{R}^{n-1},\nu)}
\leq C\, \|f\|_{L^p(\mathbb{R}^{n-1},\nu)},
\quad \text{ for every }\,f \in L^p(\mathbb{R}^{n-1},\nu).
\end{equation}
As we are assuming \eqref{eq:NT-M1}-\eqref{eq:NT-M2}
we can now use Theorem \ref{thm:XExtrapolation}
which readily gives \eqref{eq:NT}.
\end{proof}

\begin{notationconvention}\label{convention}
We adopt the following conventions from \cite{MMMM-ho} to simplify the presentation of rest of this section:
\begin{enumerate}
\item Given two vector-valued functions $f$ and $g$, we write $f\equiv g$ to indicate that each component of $f$ can be expressed as a finite linear combination of the components of $g$.
\item Given a coefficient tensor $A$, we write $M_A f$ to denote the multiplication of some (or all) components of the vector-valued function $f$ with entries from $A$ or from $(A_{nn})^{-1}$.
\item The symbol $\partial_\tau$ will denote any of the derivatives $\partial_1, ... \partial_{n-1}$, and $\partial_\tau^k$ will stand for its k-fold iteration.
\item For every $\alpha,\beta\in\{1,...,M\}$, $Q'$ will denote any $M\times M$ matrix with entries of the form $Q^{(j)}_{\alpha\beta}$ with $j \in \{1,...,n-1\}$.
\end{enumerate}
\end{notationconvention}

The proof of the following proposition can be found in \cite[Proposition 4.9]{MMMM-ho} for the case in which $f \in L_k^p(\mathbb{R}^{n-1})$ for some $p\in(1,\infty)$. A close inspection of the proof reveals that the result is actually valid whenever the convolutions $P^L_t \ast \partial_\tau ^{\,i} f$ and $Q'_t \ast \partial_\tau ^{\,i} f$ are well-defined for all $i \in \{0,1,...,k\}$. 
As a consequence of \eqref{eq:hprop-2}, this is indeed the case under the usual assumptions whenever $f\in [(\mathbb{X}_v)^w_{\ell}]^M$ and $k\in\{0,...,\ell\}$. This allows us to state the result as follows.

\begin{proposition}\label{prop:HOu}
Let $L$ be a second-order, homogeneous, elliptic, constant complex coefficient, $M\times M$ system in ${\mathbb{R}}^n$ as in \eqref{eq:L-1}-\eqref{eq:L-2} with the
property that $\mathfrak{A}_L^{\rm dis}\neq\varnothing$. Fix $\ell \in \mathbb{N}_0$, $A \in\mathfrak{A}_L^{\rm dis}$ and denote by $P^L$ the Poisson kernel for $L$ given by Theorem \ref{thm:Poisson}. Let $v$ and $w$ be weights on $(\mathbb{R}^{n-1},\mathcal{L}^{n-1})$ such that $v \in L^1_{\loc}(\mathbb{R}^{n-1})$ and $\mathbb{X}_v$ a Banach function space over $(\mathbb{R}^{n-1}, v)$. Assume that there is $C_0\in (0,\infty)$ such that 
\begin{align}
\|(\mathcal{M}h)\,w\|_{\mathbb{X}_v} &\leq C_0 \|hw\|_{\mathbb{X}_v}, \qquad \text{ for each }\,h \in \mathbb{M}, 
\\
\|(\mathcal{M}'_{v}h)\,w^{-1} \|_{\mathbb{X}'_v} &\leq C_0 \|hw^{-1}\|_{\mathbb{X}'_v}, \qquad \text{ for each }\,h \in\mathbb{M}.
\end{align}
For each $f\in [(\mathbb{X}_v)^w_{\ell}]^M$, define
\begin{equation}
u(x',t) := (P^L_t \ast f)(x'), \qquad \text{ for each } (x',t) \in \mathbb{R}_{+}^n.
\end{equation}
Then, for every $(x',t) \in \mathbb{R}_{+}^n$ and $k\in\{0,1,...,\ell\}$,
\begin{equation}\label{eq:HOu}
\nabla^k u(x',t) \equiv M_A (P^L_t \ast (M_A \partial_\tau ^k f))(x')
 + M_A (Q'_t \ast (M_A \partial_\tau^k f))(x').
\end{equation}
\end{proposition}

All necessary tools have now been established and we are ready to address the proofs of our main results.

\begin{proof}[Proof of Theorem~\ref{thm:BVP-main}]
We first prove the uniqueness of solutions. Assume that
$u$ satisfies \eqref{eq:BVP-main}. In particular, 
$(\mathcal{N}_{\kappa}u)\, w\in\mathbb{X}_v$ and hence \eqref{eq:hprop-2}
implies
\begin{equation}
\int_{\mathbb{R}^{n-1}}\frac{({\mathcal{N}}_{\kappa}u)(x')}{1+|x'|^{n-1}}\,dx'<\infty.
\end{equation}
This allows us to apply Proposition \ref{prop:uniq} to conclude that the solution is unique and it is given by \eqref{eq:BVP-main-sol}.

Next, fix a function $f \in [(\mathbb{X}_v)^w_{\ell}]^M$.
In particular, $fw \in [\mathbb{X}_v]^M$ so that \eqref{eq:hprop-2} guarantees that
\begin{equation}
\int_{\mathbb{R}^{n-1}}\frac{|f(x')|}{1+|x'|^{n-1}}\,dx'<\infty,
\end{equation}
and hence, Proposition \ref{prop:existence} 
can be applied and the function $u$ defined as in \eqref{eq:BVP-main-sol}
satisfies the first, second and last conditions in \eqref{eq:BVP-main}.

To finish the proof it remains to prove that 
$(\mathcal{N}_{\kappa}(\nabla ^k u))\, w\in\mathbb{X}_v$ for each $k\in\{0,1,...,\ell\}$
and that \eqref{eq:BVP-main-estimate} is satisfied. With this goal in mind, 
fix $k\in\{0,1,...,\ell\}$ and note that $(\partial^k_\tau f)w \in [\mathbb{X}_v]^M$.
Hence, from \eqref{eq:hprop-2} and Proposition \ref{prop:existence} we have
\begin{equation}\label{eq:existence-NuMf-inproof}
\big(\mathcal{N}_{\kappa}(P^L_t \ast \partial^k_\tau f \big))(x')\leq C\big(\mathcal{M}(\partial^k_\tau f)\big)(x'),\qquad\text{ for all }\,x'\in\mathbb{R}^{n-1}.
\end{equation}
From this and \eqref{eq:M1},
\begin{equation}\label{eq:estimate-X-inproof}
\|\big(\mathcal{N}_{\kappa}(P^L_t \ast \partial^k_\tau f)\big) w\|_{\mathbb{X}_v} \leq C
\| (\partial^k_\tau f) w\|_{[\mathbb{X}_v]^M}.
\end{equation}
Finally, using Proposition \ref{prop:HOu},
\begin{align}
\|(\mathcal{N}_{\kappa}(\nabla ^k u))\,w\|_{\mathbb{X}_v}
& \leq C \big( \| (\mathcal{N}_{\kappa} (M_A(P^L_t \ast (M_A \partial_\tau ^k f)))) \,w \|_{\mathbb{X}_v}
\nonumber \\
& \qquad \qquad + \|  (\mathcal{N}_{\kappa}(M_A(Q'_t \ast (M_A \partial_\tau^k f)))) \,w \|_{\mathbb{X}_v} \big)
\nonumber \\
& \leq C\, \| \mathcal{N}_{\kappa} (P^L_t \ast ( \partial_\tau ^k f)) \,w \|_{\mathbb{X}_v}
 + \| \mathcal{N}_{\kappa} (Q'_t \ast ( \partial_\tau ^k f)) \,w \|_{\mathbb{X}_v}
\nonumber \\
& \leq C\, \|(\partial_\tau^k f)\,w\|_{[\mathbb{X}_v]^M}
\leq C\, \|f\|_{[(\mathbb{X}_v)^w_k]^M},
\end{align}
where the third inequality follows from \eqref{eq:estimate-X-inproof} and Lemma \ref{lem:NT}.
This finishes the proof of the theorem.
\end{proof}

To finish this section, we now turn to the proof of Theorem~\ref{thm:BVP-sharp}.

\begin{proof}[Proof of Theorem~\ref{thm:BVP-sharp}]
To prove $(a)$,
fix a measurable nonnegative function $f$ so that $f w \in \mathbb{X}_v$ and
assume that the Dirichlet problem \eqref{eq:BVP-sharp1} has a 
solution $u$ satisfying \eqref{eq:BVP-sharp1-estimates}. In particular,
$(\mathcal{N}_{\kappa}u)\, w\in\mathbb{X}_v$ and hence \eqref{eq:hprop-2}
implies that 
\begin{equation}
\int_{\mathbb{R}^{n-1}}\frac{({\mathcal{N}}_{\kappa}u)(x')}{1+|x'|^{n-1}}\,dx'<\infty.
\end{equation}
Consequently, we may now apply Proposition \ref{prop:uniq} which guarantees that the solution is given by
\begin{equation}
u(x',t)=(P_t^{\Delta}\ast f)(x') \quad \text{for all }(x',t)\in{\mathbb{R}}^n_{+},
\end{equation}
where $P^{\Delta}$ is the harmonic Poisson kernel defined as in \eqref{eq:Poisson-harmonic}. Then, there exists a constant $C\in(0,\infty)$ such that
\begin{align}
u(x',t) & = (P_t^{\Delta} \ast f)(x') = \int_{\mathbb{R}^{n-1}} t^{1-n} P^{\Delta} \Big( \frac{x'-y'}{t} \Big) f(y')\,dy'
\nonumber\\[4pt]
&  = \frac{2}{\omega_{n-1}}  \int_{\mathbb{R}^{n-1}} \frac{t}{(t^2 + |x'-y'|^2)^{n/2}} f(y') \, dy'
\nonumber\\[4pt]
& \geq C \int_{B(x',t)} \frac{t}{(t^2 + t^2)^{n/2}} f(y') \, dy'
\nonumber\\[4pt]
& = C \fint_{B(x',t)} f(y') \, dy',
\end{align}
where $B(x',t)$ denotes the ball centered at $x'$ with radius $t$ in $\mathbb{R}^{n-1}$. In turn, this implies
\begin{equation}
(\mathcal{M} f)(x') \leq C \,
\sup_{t>0} u(x',t) \leq C
(\mathcal{N}_{\kappa}u)(x')
\qquad \text{ for each }x' \in \mathbb{R}^{n-1}.
\end{equation}
As a consequence, thanks to \eqref{eq:BVP-sharp1-estimates}, we obtain
\begin{equation}
\| (\mathcal{M} f) \,w \|_{\mathbb{X}_v}
\leq C \| (\mathcal{N}_{\kappa}u)\,w \|_{\mathbb{X}_v}
\leq C \| f \,w \|_{\mathbb{X}_v},
\end{equation}
that is, \eqref{eq:sharp1-M1} holds.

For the proof of $(b)$, we proceed analogously. Fix a measurable function $f$ so that $f w^{-1}v^{-1} \in \mathbb{X}_v$ and,
reasoning as before, the well-posedness of the Dirichlet problem 
\eqref{eq:BVP-sharp2} satisfying \eqref{eq:BVP-sharp2-estimates} implies 
\begin{equation}
\| (\mathcal{M} f) \,w^{-1}v^{-1} \|_{\mathbb{X}'_v}
\leq C \| (\mathcal{N}_{\kappa}u)\,w^{-1}v^{-1} \|_{\mathbb{X}'_v}
\leq C \| f \,w^{-1}v^{-1} \|_{\mathbb{X}'_v}.
\end{equation}
By replacing $f$ by $hv$, one can check that
this is actually equivalent to the fact that \eqref{eq:sharp2-M2} holds for
every measurable function $h$ so that $hw^{-1} \in \mathbb{X}_v$.

Finally, $(c)$ is a direct consequence of $(a)$, $(b)$, Theorem \ref{thm:BVP-main}, and Remark \ref{rem:dual}.
\end{proof}

\section{Particular cases}\label{sec:examples}

In order to illustrate the breadth of Theorem \ref{thm:BVP-main}, we consider here a variety of function spaces for which the theorem is valid and establish the corresponding new well-posedness results generated in these settings.

\subsection{Weighted Lebesgue spaces}\label{subsec:exa-Lp}

The most basic scenario is the case of Lebesgue spaces.
Specifically, fix $p\in (1,\infty)$ and $w\in A_p(\mathbb{R}^{n-1})$ a Muckenhoupt weight.
We will apply Theorem \ref{thm:BVP-main} with $v\equiv 1$, $\mathbb{X}_v = L^p(\mathbb{R}^{n-1})$ and $w^{1/p}$ in place
of $w$.
In this case,
\begin{equation}
\mathbb{X}'_v=L^{p'}(\mathbb{R}^{n-1}), \quad \text{where} \quad
1/p+1/p'=1,
\end{equation}
and
\begin{equation}
(\mathbb{X}_v)^{w^{1/p}}_\ell = L_\ell^{p}(\mathbb{R}^{n-1},w) =
\{ f \in L^{p}(\mathbb{R}^{n-1},w) : \partial^{\alpha} f \in L^{p}(\mathbb{R}^{n-1},w)
\text{ if } |\alpha| \leq \ell \},
\end{equation}
is the usual weighted $L^p$-based Sobolev space.
Moreover, in this scenario, the hypothesis \eqref{eq:M1}-\eqref{eq:M2} amounts to saying that the Hardy-Littlewood maximal function $\mathcal{M}$ is bounded both on $L^{p}(\mathbb{R}^{n-1},w)$
and on $L^{p'}(\mathbb{R}^{n-1},w^{-p'/p})=L^{p'}(\mathbb{R}^{n-1},w^{1-p'})$.
As a consequence of \eqref{eq:Ap-p1} and \eqref{eq:Mbounded-Lpw},
both of these conditions are actually equivalent to the fact that 
$w\in A_p(\mathbb{R}^{n-1})$.
This leads to the following well-posedness result.

\begin{example}\label{ex:Lpw}
Let $L$ be a second-order, homogeneous, elliptic, constant complex coefficient, $M\times M$ system in ${\mathbb{R}}^n$ as in \eqref{eq:L-1}-\eqref{eq:L-2}
such that $\mathfrak{A}_L^{\rm dis} \neq \varnothing$. Fix an aperture parameter $\kappa>0$ and $\ell \in \mathbb{N}_0$. Then, for every $p\in(1,\infty)$ and $w \in A_p(\mathbb{R}^{n-1})$, the $\ell$-th order regularity problem
\begin{equation}\label{eq:BVP-Lpw}
\left\lbrace
\begin{array}{l}
u\in\big[\mathscr{C}^{\infty}(\mathbb{R}^n_{+})\big]^M,
\\[4pt]
Lu=0\,\,\text{ in }\,\,\mathbb{R}^n_{+},
\\[4pt] 
(\mathcal{N}_{\kappa}(\nabla ^k u)) \in L^p(\mathbb{R}^{n-1},w), \quad k\in\{0,1,...,\ell\},
\\[4pt]
\restr{u}{\partial\mathbb{R}^n_{+}}^{{}^{\kappa-{\rm n.t.}}}=f \in \big[L_\ell^p(\mathbb{R}^{n-1},w)\big]^M,
\end{array}
\right.
\end{equation}
has a unique solution given by $u(x',t)=(P_t^L*f)(x')$ for every $(x',t)\in{\mathbb{R}}^n_{+}$, where $P^L$ denotes the Poisson kernel for $L$ in $\mathbb{R}^n_+$ from Theorem~\ref{thm:Poisson}. Furthermore, there is a constant $C\in[1,\infty)$ such that 
\begin{equation}\label{eq:BVP-Lpw-estimate}
\sum_{k=0}^{\ell}\|(\mathcal{N}_{\kappa}(\nabla ^k u))\|_{L^p(\mathbb{R}^{n-1},w)} \leq C\, \|f\|_{[L_\ell^p(\mathbb{R}^{n-1},w)]^M}.
\end{equation}
\end{example}

Note that if we set $w\equiv 1$ then the previous example becomes the $\ell$-th order regularity problem for ordinary (unweighted) Lebesgue spaces, which was studied in \cite{MMMM-ho}.

Finally, it is worth emphasizing that this example illustrates the necessity of introducing a second weight $w$ in the formulation of Theorem \ref{thm:BVP-main}. Indeed, although $\mathbb{X}_v=L^p(\mathbb{R}^{n-1},v)$ with $v \in A_p(\mathbb{R}^{n-1})$ is a Banach function space over $(\mathbb{R}^{n-1},v)$, in this setting the associate space must be computed with respect to the measure space $(\mathbb{R}^{n-1},v)$, yielding $\mathbb{X}'_v=L^{p'}(\mathbb{R}^{n-1},v)$. However, the Hardy-Littlewood maximal operator $\mathcal{M}$ may fail to be bounded on this space and, consequently, hypothesis \eqref{eq:M1} is not guaranteed to hold.

\subsection{Weighted rearrangement invariant Banach function spaces}\label{exa-RIBFS}

Let $\mathcal{L}_f$ denote the distribution function of $f \in \mathbb{M}$, that is, 
\begin{equation}
\mathcal{L}_f (\lambda) := \mathcal{L}^{n-1} ( \{ x' \in \mathbb{R}^{n-1} : |f(x')|>\lambda \} ) \text{,     } \qquad \lambda \geq 0. 
\end{equation}
Let $\mathbb{X}$ be a Banach function space over $(\mathbb{R}^{n-1},\mathcal{L}^{n-1})$, that is, a space as in Definition \ref{def:BFS} with $v \equiv 1$, a case in which we agree to omit the subscript and write $\mathbb{X}$ in place of $\mathbb{X}_v$. The space $\mathbb{X}$
is called {\tt rearrangement invariant} if $\norm{ f } _{\mathbb{X}} = \norm{ g } _{\mathbb{X}}$ 
for every pair of functions $f,g \in \mathbb{X}$ such that $\mathcal{L}_f = \mathcal{L}_g$.

Given $f \in \mathbb{M}$, the decreasing rearrangement of $f$ with 
respect to $\mathcal{L}^{n-1}$ is the function defined as 
\begin{equation}
f^{*}(t) := \inf \{ \lambda \geq 0 : \mathcal{L}_f (\lambda) \leq t \}, \qquad t\in[0,\infty).
\end{equation} 
Note that $f$ and $f^{*}$ have the same distribution function. This leads to the
Luxemburg representation theorem: given a rearrangement invariant Banach function space $\mathbb{X}$, there exists another rearrangement invariant Banach function space
$\overline{\mathbb{X}}$ on $[0,\infty)$ such that $f \in \mathbb{X}$ if and only if 
$f^{*}\in \overline{\mathbb{X}}$ and 
$\norm{ f^{*} } _{\overline{\mathbb{X}}} = \norm{ f } _{\mathbb{X}}$. 

Next, we define the lower and upper Boyd indices as
\begin{equation}
p_{\mathbb{X}} := \lim_{t \to \infty} 
\frac{\log t}{\log \norm{ D_t } _{\overline{\mathbb{X}}\to \overline{\mathbb{X}}}} \text{,     } 
\qquad q_{\mathbb{X}} 
:= \lim_{t \to 0^{+}} \frac{\log t}{\log \norm{ D_t } _{\overline{\mathbb{X}}\to \overline{\mathbb{X}}}}, 
\end{equation}
where $D_t: \overline{\mathbb{X}} \to \overline{\mathbb{X}}$ is the dilation operator defined by $D_t f(s):=f(s/t)$ for every $s \geq 0$ and $t>0$.
Boyd indices capture essential information about the interpolation of operators in rearrangement invariant Banach function spaces (cf. \cite{BS88}).

Our next goal is to mirror the definition of rearrangement invariant spaces
to introduce a weighted version of them. Specifically, given a weight $v$ and $f\in\mathbb{M}$, the distribution function of $f$ with respect to $v$ is defined as 
\begin{equation}
v_f (\lambda) := v( \{ x' \in \mathbb{R}^{n-1} : |f(x')|>\lambda \} ) \text{,     } \qquad \lambda \geq 0,
\end{equation}
and the decreasing rearrangement of $f$ with respect to $v$ is
\begin{equation}
f_v^{*}(t) := \inf \{ \lambda \geq 0 : v_f (\lambda) \leq t \}, \qquad t\in[0,\infty).
\end{equation} 
This provides a way to define weighted rearrangement invariant Banach function spaces as
\begin{equation} 
\mathbb{X}(v) := \{ f \in \mathbb{M} : \norm{ f^{*}_{v} } _{\overline{\mathbb{X}}} < \infty \}, 
\end{equation}
which is a Banach function space on $(\mathbb{R}^{n-1},v)$ with norm 
$\norm{ f } _{\mathbb{X}(v)} = \norm{ f^{*}_{v} } _{\overline{\mathbb{X}}}$.
As a particular case, note that if $\mathbb{X}=L^p(\mathbb{R}^{n-1})$ for $p \in (1,\infty)$ then $\mathbb{X}(v)=L^p(\mathbb{R}^{n-1},v)$.

We can now apply Theorem \ref{thm:BVP-main} for $\mathbb{X}_v=\mathbb{X}(v)$. In this setting,
\begin{equation}
(\mathbb{X}_v)^{w}_\ell =  (\mathbb{X}(v))^{w}_\ell = 
\{ f \in \mathbb{M} : f w \in \mathbb{X}(v) \text{ and } (\partial^{\alpha} f) w \in \mathbb{X}(v)
\text{ if } |\alpha| \leq \ell \}
\end{equation}
with the natural associated norm.
For the hypothesis \eqref{eq:M1}-\eqref{eq:M2} to hold in this context, we need to impose some conditions on the weights. It was proved in \cite{CMM22} that whenever
$1<p_{\mathbb{X}}\leq q_{\mathbb{X}} < \infty$ we have
\begin{equation}
\parbox{13cm}{\begin{center}
$\norm{(\mathcal{M}h)\, w}_{\mathbb{X}(v)} \leq C_0\, \|hw\|_{\mathbb{X}(v)}$ and $\|(\mathcal{M}'_{v} h)\, w^{-1} \|_{\mathbb{X}'(v)} \leq C_0\, \|hw^{-1}\|_{\mathbb{X}'(v)}$
whenever
$w^{p_{\mathbb{X}}}v \in A_{p_{\mathbb{X}}}(\mathbb{R}^{n-1})$, 
$w^{q_{\mathbb{X}}}v \in A_{q_{\mathbb{X}}}(\mathbb{R}^{n-1})$,  and 
$v \in A_{\infty}(\mathbb{R}^{n-1})$.
\end{center}}
\end{equation} 
This leads to the following well-posedness result.

\begin{example}
Let $L$ be a second-order, homogeneous, elliptic, constant complex coefficient, $M\times M$ system in ${\mathbb{R}}^n$ as in \eqref{eq:L-1}-\eqref{eq:L-2} such that $\mathfrak{A}_L^{\rm dis} \neq \varnothing$. Fix an aperture parameter $\kappa>0$ and $\ell \in \mathbb{N}_0$. Let $\mathbb{X}$ be a rearrangement
invariant Banach function space over $(\mathbb{R}^{n-1},\mathcal{L}^{n-1})$
such that $1<p_{\mathbb{X}}\leq q_{\mathbb{X}} < \infty$ and let $v,w$ be
a pair of weights such that 
\begin{equation}\label{eq:BVP-ri-weights}
w^{p_{\mathbb{X}}}v \in A_{p_{\mathbb{X}}}(\mathbb{R}^{n-1}), \quad
w^{q_{\mathbb{X}}}v \in A_{q_{\mathbb{X}}}(\mathbb{R}^{n-1}),
\quad  \text{ and } \quad 
v \in A_{\infty}(\mathbb{R}^{n-1}).
\end{equation}
Then, the $\ell$-th order regularity problem
\begin{equation}\label{eq:BVP-ri}
\left\lbrace
\begin{array}{l}
u\in\big[\mathscr{C}^{\infty}(\mathbb{R}^n_{+})\big]^M,
\\[4pt]
Lu=0\,\,\text{ in }\,\,\mathbb{R}^n_{+},
\\[4pt] 
(\mathcal{N}_{\kappa}(\nabla ^k u))\, w\in\mathbb{X}(v), \quad k\in\{0,1,...,\ell\},
\\[4pt]
\restr{u}{\partial\mathbb{R}^n_{+}}^{{}^{\kappa-{\rm n.t.}}}=f \in \big[(\mathbb{X}(v))^w_\ell \big]^M,
\end{array}
\right.
\end{equation}
has a unique solution given by $u(x',t)=(P_t^L*f)(x')$ for every $(x',t)\in{\mathbb{R}}^n_{+}$, where $P^L$ denotes the Poisson kernel for $L$ in $\mathbb{R}^n_+$ from Theorem~\ref{thm:Poisson}. Furthermore, there is a constant $C\in[1,\infty)$ such that 
\begin{equation}\label{eq:BVP-ri-estimate}
\sum_{k=0}^{\ell}\|(\mathcal{N}_{\kappa}(\nabla ^k u))\,w\|_{\mathbb{X}(v)} \leq C\, \|f\|_{[(\mathbb{X}(v))^w_\ell]^M}.
\end{equation}
\end{example}

We first emphasize that if we set $w\equiv 1$ then \eqref{eq:BVP-ri} is
the $\ell$-th order regularity problem for $\mathbb{X}(v)$ and the hypothesis
\eqref{eq:BVP-ri-weights} becomes simply $v \in A_{p_{\mathbb{X}}}(\mathbb{R}^{n-1})$.

Many function spaces fall within the scope of this example.
Apart from ordinary and weighted Lebesgue spaces, 
this covers the cases 
\begin{equation}
\mathbb{X}=L^p(\mathbb{R}^{n-1}) + L^q(\mathbb{R}^{n-1})
\quad\text{ and }\quad\mathbb{X}=L^p(\mathbb{R}^{n-1}) \cap L^q(\mathbb{R}^{n-1})
\end{equation}
with $p,q\in(1,\infty)$, for both of which
$p_{\mathbb{X}} = \min\{p,q\}$ and $q_{\mathbb{X}}=\max\{p,q\}$, as well as their weighted counterparts if $w^p v \in A_p(\mathbb{R}^{n-1})$, $w^q v \in A_q(\mathbb{R}^{n-1})$ and $v \in A_{\infty}(\mathbb{R}^{n-1})$.
It also includes other Orlicz spaces such as Zygmund spaces 
\begin{equation}
\mathbb{X}=L^{p}(\log L)^{\alpha}(\mathbb{R}^{n-1})
\end{equation}
with $p\in(1,\infty)$ and $\alpha \in \mathbb{R}$, for which
$p_{\mathbb{X}} = q_{\mathbb{X}}=p$, as well as its weighted counterpart if $w^p v \in A_p(\mathbb{R}^{n-1})$ and $v \in A_{\infty}(\mathbb{R}^{n-1})$. 
Further details concerning these spaces and their Boyd indices may be found in \cite[Chapter 4]{CMP11}.

\subsection{Weighted variable exponent Lebesgue spaces}\label{subsec:ex-Lpvar}

Let $p(\cdot): \mathbb{R}^{n-1} \to (1, \infty)$ be a Lebesgue measurable function.  
The variable exponent Lebesgue space $L^{p(\cdot)}(\mathbb{R}^{n-1})$ is the Banach function space associated with the norm
\begin{equation}\label{eq:Lpvar-norm}
\|f\|_{L^{p(\cdot)}(\mathbb{R}^{n-1})} 
:= \inf \left\lbrace \lambda >0: \,\, \int_{\mathbb{R}^{n-1}}
\left( \frac{|f(x')|}{\lambda} \right)^{p(x')}
dx' \leq 1 \right\rbrace .
\end{equation}
This class of spaces extends the family of the usual Lebesgue spaces.  
Indeed, if $p(x') = p_0$ is constant, then
$L^{p(\cdot)}(\mathbb{R}^{n-1}) = L^{p_0}(\mathbb{R}^{n-1})$.
Throughout, we adopt the notation $p(\cdot)$ instead of simply $p$ to emphasize that the exponent is variable. For a given variable exponent $p(\cdot)$, we define the essential infimum and supremum as
\begin{equation}
p_{-} := \operatorname*{ess\,inf}_{x'\in\mathbb{R}^{n-1}} p(x'),
\qquad
p_{+} := \operatorname*{ess\,sup}_{x'\in\mathbb{R}^{n-1}} p(x').
\end{equation}
We say that $p(\cdot) \in {\rm LH}(\mathbb{R}^{n-1})$ if there are constants $C \in [0,\infty)$ and $p_{\infty} \in [0,\infty)$ such that for all $x', y' \in \mathbb{R}^{n-1}$, the following inequalities hold:
\begin{equation}
|p(x') - p(y')| \leq \frac{C}{-\log |x' - y'|}, \quad \text{whenever } 0 < |x' - y'| \leq \tfrac{1}{2},
\end{equation}
and
\begin{equation}
|p(x') - p_{\infty}| \leq \frac{C}{\log (e + |x'|)}.
\end{equation}
Next, we define the class of $A_{p(\cdot)}(\mathbb{R}^{n-1})$ weights as the family of weights satisfying
\begin{equation}\label{eq:Apvar}
[w]_{A_{p(\cdot)}}:=\sup_{B} \big( [\mathcal{L}^{n-1}(B)]^{-1}\,
\norm{\mathbf{1}_B w}_{L^{p(\cdot)}(\mathbb{R}^{n-1})}
\|\mathbf{1}_B w^{-1}\|_{L^{p'(\cdot)}(\mathbb{R}^{n-1})}\big)<\infty,
\end{equation}
where the supremum runs over all balls $B$ in $\mathbb{R}^{n-1}$, $\mathbf{1}_B$ denotes the characteristic function of the ball $B$, and $1/p(\cdot)+1/p'(\cdot)=1$.
Next, we recall that weighted variable exponent Lebesgue spaces $L^{p(\cdot)}(\mathbb{R}^{n-1},w)$ are defined through the norm
\begin{equation}
\| f \| _{L^{p(\cdot)}(\mathbb{R}^{n-1},w)} :=
\| f w \| _{L^{p(\cdot)}(\mathbb{R}^{n-1})}.
\end{equation}
Under the background assumption that $1<p_{-}\leq p_{+}<\infty$ and $p(\cdot)\in{\rm LH}(\mathbb{R}^{n-1})$, it is proved in \cite{CDH11, CFN12} that
\begin{align}\label{eq:Apvar-equiv}
w \in A_{p(\cdot)}(\mathbb{R}^{n-1}) \,\, &\Longleftrightarrow \,\, \|Mh\|_{L^{p(\cdot)}(\mathbb{R}^{n-1},w)} \leq C \|h\|_{L^{p(\cdot)}(\mathbb{R}^{n-1},w)}, 
\nonumber\\
 &\Longleftrightarrow \,\, \|Mh\|_{L^{p'(\cdot)}(\mathbb{R}^{n-1},w^{-1})} \leq C \|h\|_{L^{p'(\cdot)}(\mathbb{R}^{n-1},w^{-1})}.
\end{align}

With these equivalences in mind, we will now apply Theorem \ref{thm:BVP-main} for $v\equiv 1$ and $\mathbb{X}_v = L^{p(\cdot)}(\mathbb{R}^{n-1})$. In this context,
\begin{equation}
\mathbb{X}'_v=L^{p'(\cdot)}(\mathbb{R}^{n-1}), \quad \text{where} \quad
1/p(\cdot)+1/p'(\cdot)=1,
\end{equation}
and
\begin{align}
(\mathbb{X}_v)^{w}_\ell & = L_\ell^{p(\cdot)}(\mathbb{R}^{n-1},w) 
\nonumber\\
& =
\{ f \in \mathbb{M} : f \in L^{p(\cdot)}(\mathbb{R}^{n-1},w) 
 \text{ and } \partial^{\alpha} f \in L^{p(\cdot)}(\mathbb{R}^{n-1},w) 
\text{ if } |\alpha| \leq \ell \}
\end{align}
with the natural associated norm. We obtain the following result.

\begin{example}\label{ex:Lpvarw}
Let $L$ be a second-order, homogeneous, elliptic, constant complex coefficient, $M\times M$ system in ${\mathbb{R}}^n$ as in \eqref{eq:L-1}-\eqref{eq:L-2} such that $\mathfrak{A}_L^{\rm dis} \neq \varnothing$. Fix an aperture parameter $\kappa>0$ and $\ell \in \mathbb{N}_0$. Assume that $p(\cdot)$ is a variable exponent such that $1 < p_{-} \leq p_{+} < \infty$, $p(\cdot) \in {\rm LH}(\mathbb{R}^{n-1})$ and $w \in A_{p(\cdot)}(\mathbb{R}^{n-1})$.
Then, the $\ell$-th order regularity problem
\begin{equation}\label{eq:BVP-Lpvarw}
\left\lbrace
\begin{array}{l}
u\in\big[\mathscr{C}^{\infty}(\mathbb{R}^n_{+})\big]^M,
\\[4pt]
Lu=0\,\,\text{ in }\,\,\mathbb{R}^n_{+},
\\[4pt] 
\mathcal{N}_{\kappa}(\nabla ^k u) \in L^{p(\cdot)}(\mathbb{R}^{n-1},w), \quad k\in\{0,1,...,\ell\},
\\[4pt]
\restr{u}{\partial\mathbb{R}^n_{+}}^{{}^{\kappa-{\rm n.t.}}}=f \in \big[ L_\ell^{p(\cdot)}(\mathbb{R}^{n-1},w) \big]^M,
\end{array}
\right.
\end{equation}
has a unique solution given by $u(x',t)=(P_t^L*f)(x')$ for every $(x',t)\in{\mathbb{R}}^n_{+}$, where $P^L$ denotes the Poisson kernel for $L$ in $\mathbb{R}^n_+$ from Theorem~\ref{thm:Poisson}. Furthermore, there is a constant $C\in[1,\infty)$ such that 
\begin{equation}\label{eq:BVP-Lpvarw-estimate}
\sum_{k=0}^{\ell}\|\mathcal{N}_{\kappa}(\nabla ^k u)\|_{L^{p(\cdot)}(\mathbb{R}^{n-1},w)} \leq C\, \|f\|_{[L_\ell^{p(\cdot)}(\mathbb{R}^{n-1},w)]^M}.
\end{equation}
\end{example}

As for the particular case $w\equiv 1$, one can see from \cite[Theorem 2.4]{CW17} that
\begin{align}
w \equiv 1 \in A_{p(\cdot)}(\mathbb{R}^{n-1}) \,\, \text{ for every } \,\,p(\cdot) \in{\rm LH}(\mathbb{R}^{n-1}) \text{ with } 1<p_{-}\leq p_{+}<\infty. 
\end{align}
As a consequence, we obtain well-posedness for the $\ell$-th order regularity problem \eqref{eq:BVP-Lpvarw} with $w\equiv 1$.


\def\ocirc#1{\ifmmode\setbox0=\hbox{$#1$}\dimen0=\ht0 \advance\dimen0
by1pt\rlap{\hbox to\wd0{\hss\raise\dimen0
\hbox{\hskip.2em$\scriptscriptstyle\circ$}\hss}}#1\else {\accent"17 #1}\fi}
\def\cprime{$'$} \def\cprime{$'$} \def\cprime{$'$}
\end{document}